\documentclass[11pt]{amsart}
\usepackage{amsmath, amsthm, amssymb}
\usepackage{graphicx}
\usepackage{mathbbol}
\usepackage{txfonts}
\usepackage[usenames]{color}
\usepackage{srcltx} % JD - to jump between source and preview
\newtheorem{thm}{Theorem}[section]

\newtheorem{ex}[thm]{Example}
\newtheorem{lem}[thm]{Lemma}   %remember to switch all {lem} to {lem}
\newtheorem{prop}[thm]{Proposition}
\newtheorem{defn}[thm]{Definition}
\newtheorem{rmrk}[thm]{Remark}   %remember to switch all {rmrk} to {rmrk}

\newcommand{\be}{\begin{equation}}
\newcommand{\ee}{\end{equation}}

\newcommand{\ptGHto}{\stackrel { \textrm{ptGH}}{\longrightarrow} }
\newcommand{\GHto}{\stackrel { \textrm{GH}}{\longrightarrow} }
\newcommand{\Fto}{\stackrel {\mathcal{F}}{\longrightarrow} }
\newcommand{\vare}{\varepsilon}

\newcommand{\diam}{\operatorname{Diam}}

\newcommand{\Lip}{\operatorname{Lip}}

\newcommand{\mass}{{\mathbf M}}

\newcommand{\dil}{\textrm{dil}}

\newcommand{\vol}{\operatorname{Vol}}

\begin{document}

\title{SWIF Convergence of Smocked Metric Spaces}

\author{M. Dinowitz}
\address{M. Dinowitz (CUNYGC)}
\author{H. Drillick}
\address{H. Drillick (Columbia)}
\author{M. Farahzad}
\address{M. Farahzad (U Toronto)}
\author{C. Sormani}\thanks{The research was funded in part by NSF DMS-1612049 }
\address{C. Sormani (Lehman and CUNYGC)}
\author{A. Yamin}
\address{A. Yamin (Stony Brook)}


%49Q15 (Geometric measure and integration theory, integral and normal currents)

%\subjclass[2000]{49Q15}

\begin{abstract} 
In this paper we explore a special class
of metric spaces called smocked metric spaces and study their tangent cones at infinity.
We prove that under the right hypotheses, the rescaled limits of balls converge in both the 
Gromov-Hausdorff and Intrinsic Flat sense to normed spaces.     This paper will
be applied in upcoming work by Kazaras and Sormani concerning
Gromov's conjectures on the properties of GH and SWIF limits of 
Riemannian manifolds with positive scalar curvature.
\end{abstract}

\maketitle

\section{Introduction}

In 1983 Gromov introduced the notion Gromov-Hausdorff (GH) convergence of Riemannian manifolds and metric spaces \cite{Gromov-metric}.   Gromov proved all GH limits are geodesic metric spaces but they may not have the same dimension as the sequence. In \cite{SW-JDG} Sormani and Wenger introduced the notion of intrinsic flat (SWIF) convergence of Riemannian manifolds. SWIF limit spaces are called integral current spaces: they are countably $\mathcal{H}^m$ rectifiable metric spaces with normed tangent spaces almost everywhere of the same dimension as the original sequence.   When the SWIF and GH limits agree, then GH limits have far more structure than initially proven by Gromov.  

Here we prove that balls in smocked metric spaces are integral current spaces [Theorem~\ref{T-B-R}], and we explore the SWIF limits of these balls under rescaling.  We prove that under certain hypotheses the rescaled balls converge in both the GH and the SWIF sense to the same limit space and that this space is a normed metric space [Theorem~\ref{thm-SWIF=GH-R}].   
In future work of Kazaras and Sormani \cite{Kazaras-Sormani-tori},  the theorems here will be applied to 
address questions regarding the SWIF and GH limits of Riemannian manifolds with scalar curvature bounds (cf.\cite{Gromov-Sormani-IAS-report}\cite{Sormani-Scalar-21}).

\begin{figure}[h]
\includegraphics[width=.3 \textwidth]{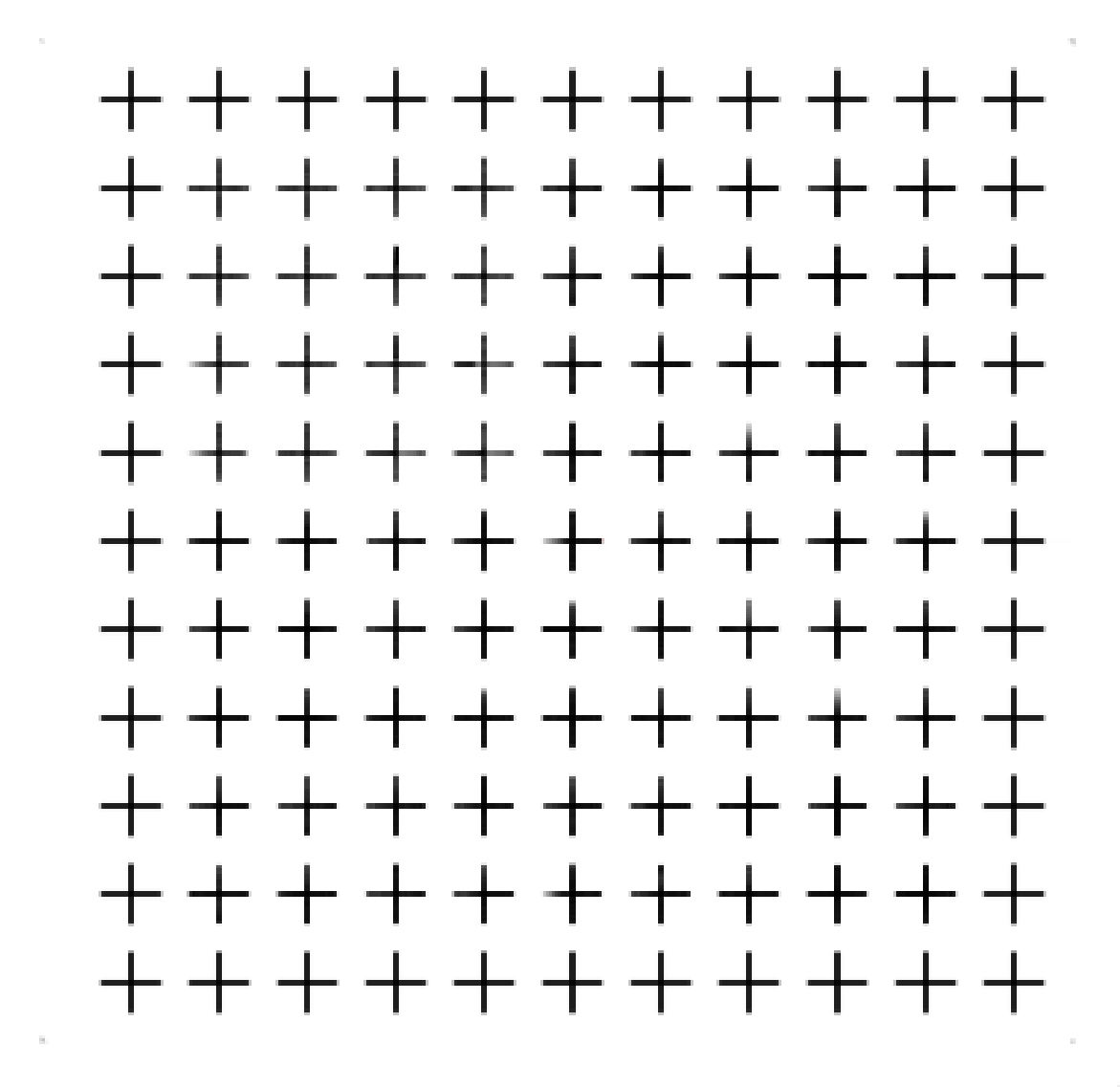} 
\includegraphics[width=.3 \textwidth]{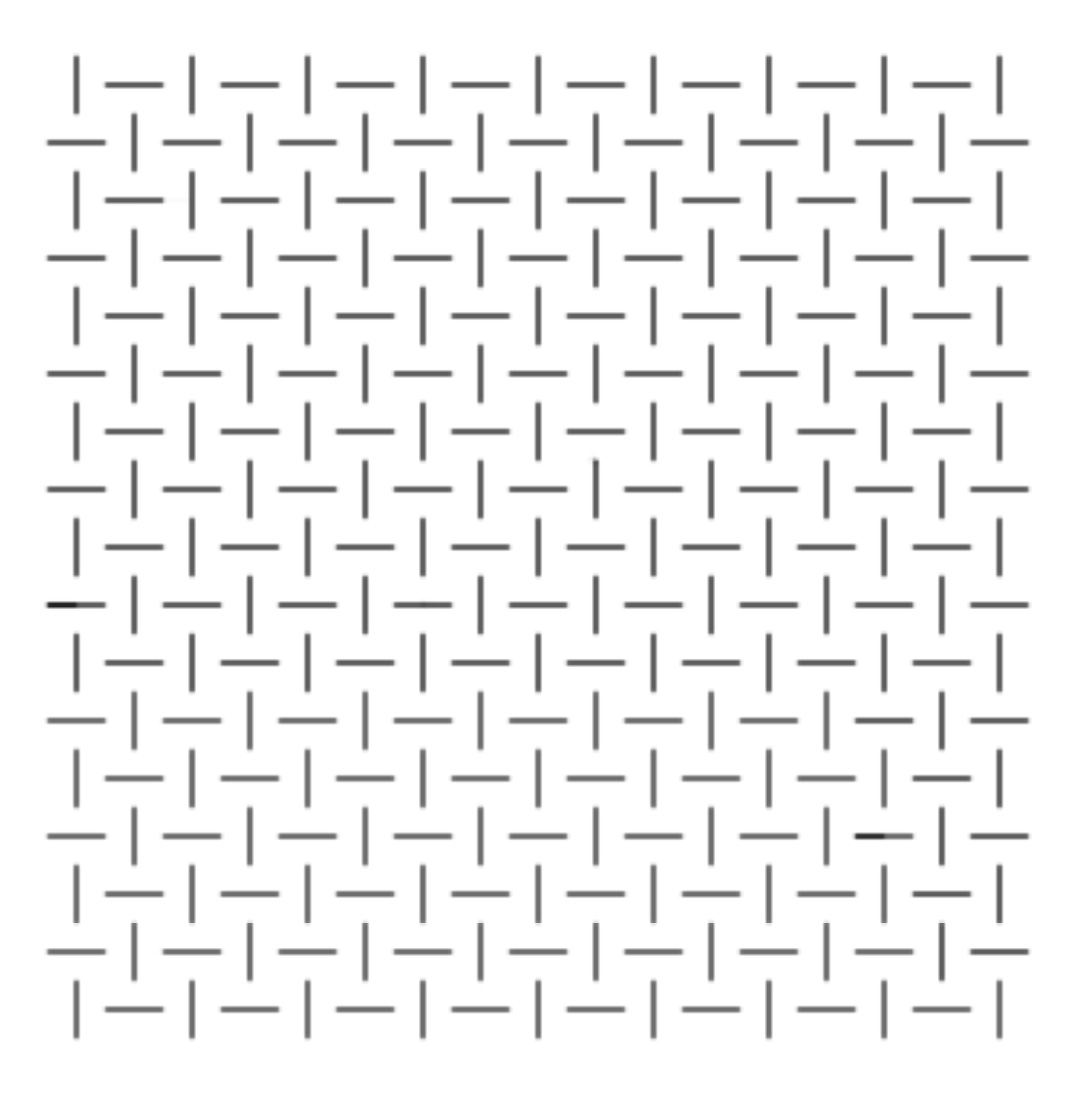} 
\includegraphics[width=.3 \textwidth]{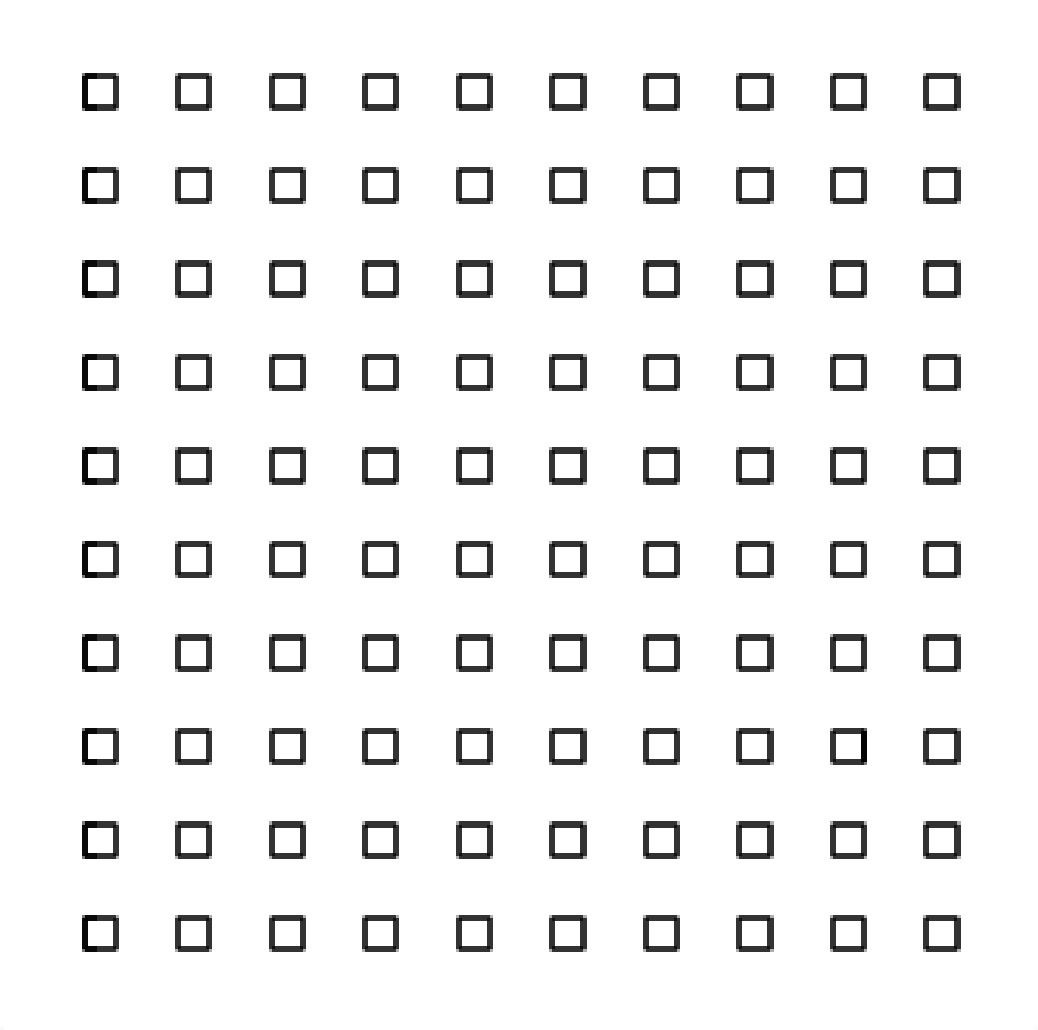} \\
\caption{Smocking patterns here were used to define the nice smocked metric spaces: $X_+$, $X_T$, and $X_\square$
in \cite{Smocked}.}
\label{fig:smocking-patterns}
\end{figure}

Recall that a smocked metric space is a metric space created by selecting a collection (called a smocking pattern) of disjoint compact subsets (called stitches) in Euclidean space and identifying all the points in each stitch to a single point (cf. Definition~\ref{defn-smock}within).   See Figure~\ref{fig:smocking-patterns} for a variety of smocking patterns of stitches on a Euclidean plane that have been analyzed by the authors and their collaborators in \cite{Smocked}.   The patterns of the stitches need not be periodic but we do require that stitches be separated from one another.   The stitches do not have to be one dimensional but we do say the space is ``nice'' if the volumes of the tubular neighborhoods of the patterns behave well (cf. Definition~\ref{smocking-nice}).    

Note that any smocked metric space, $(X,d)$, can be identified with a psuedometric on Euclidean space $({\mathbb E}^N, \bar{d})$ where $\bar{d}(p,q)=0$ whenever $p,q$ lie in the same stitch.  The distances between points $p,q$ that do not lie on the same stitch are found by taking the minimum length of straight line segments running between them jumping across the stitches.  We review the definition of these spaces and properties of these spaces in Section~\ref{ReviewSmocked}.  It is worth noting that these spaces are easy for even undergraduates to understand and there are a variety of undergraduate and masters level research projects suggested at the end of \cite{Smocked}.   

The notion of an integral current space is far more difficult to understand as it relies on the theory of Ambrosio-Kirchheim developed in \cite{AK}.   An integral current space is an integer rectifiable metric space (which is a space covered by a countable collection of biLipschitz charts with integer weights) that has a boundary which is also integer rectifiable \cite{SW-JDG}.   
In Section~\ref{CurrentSpaces}, we review this theory while proving that balls within smocked metric spaces are integral current spaces.  See Theorem~\ref{T-B-R}.

The SWIF distance between a pair of integral current spaces, $M_1$ and $M_2$, is estimated by finding a higher dimensional integral current space, $Z$, and distance preserving maps $\varphi_i: M_i \to Z$, and then estimating the volumes between their images.   In Section~\ref{SWIF}, we review this theory while estimating the SWIF distance between balls in smocked metric spaces balls in normed spaces.    We prove that under certain hypotheses on the pseudometrics, the rescaled limits of the balls in a smocked metric spaces converge in both the GH and SWIF sense to a unique tangent cone at infinity that is a normed space:

\begin{thm}\label{thm-SWIF=GH-R}
Suppose we have a nice  smocked metric space, $(X, d)$, 
 such that
\be
|\,\bar{d}(x, x')\,- \,[F(x)-F(x')] \,| \,\le \, K \qquad \forall x,x' \in {\mathbb{E}}^N
\ee
where $F: {\mathbb{E}}^N \to [0,\infty)$ is a norm.

Then $(X,d)$ has a unique GH=SWIF tangent cone at infinity, 
$({\mathbb{R}}^N, d_F)$,
where
\be
d_F(x,x')=||x-x'||_F=F(x-x').
\ee
That is, for any basepoint $x_0\in X$, $r,R>0$
the balls $\bar{B}(x_0,Rr)$
viewed as integral current spaces as in 
rescaled by $R$
converge in both the GH and SWIF sense to a ball 
of radius $r$ in the normed space $({\mathbb{R}}^N, d_F)$.
\end{thm}

We hope that this paper will be understandable even to the student and that it can be used to provide an introduction
to both Gromov-Hausdorff and Intrinsic Flat Convergence of metric spaces.  We close the paper by applying the main theorem 
to the smocked metric spaces that were analyzed by the authors, Kazaras, Afrifa, Antonetti, George, Hepburn, Huynh, Minichiello, and Rendla in \cite{Smocked}, \cite{Smocked-H}, and \cite{Smocked-Equals}.

\section{Background on Smocked Spaces} \label{ReviewSmocked}

\subsection{Review of Smocked Spaces}

The notion of a smocked space was recently introduced by Sormani, Kazaras, and their team of students
in \cite{Smocked}.  The notion is built upon the classical handcraft of smocked fabrics.  One has a pattern 
of intervals on a cloth and each interval is sewn with a thread and pulled to a point.  A smocked space is
similarly a plane of some dimension with a collection of intervals each of which is identified to a point.
See Figure~\ref{fig:smocking-patterns} for examples of patterns studied in \cite{Smocked} and
here.

\begin{defn} \label{defn-smock}
Given a Euclidean space, $\mathbb{E}^N$, and a finite or countable collection of
disjoint connected compact sets called {\bf smocking intervals}, 
\be
\mathcal{I}=\{I_j: \, j \in J\},
\ee
with a positive {\bf smocking separation factor},
\be \label{s-factor}
\delta=\min\{|z-z'|: \, z\in I_j, \, z'\in I_{j'},\, j\neq j' \in J\} >0,
\ee
 one can define the {\bf smocked
metric space}, $(X,d)$, in which each interval is viewed as a single point.  
\be
X = \left\{ x: \, x \in {\mathbb{E}^N}\setminus S\right\} \cup \mathcal{I}
\ee
where $S$ is the {\bf smocking set} or {\bf smocking pattern}:
\be
S= \bigcup_{j \in J} I_j .
\ee
We have a {\bf smocking map} 
$
\pi: \mathbb{E}^N \to X $ defined by
\be
 \pi(x) = \begin{cases} 
          x & \textrm{ for }x \in {\mathbb{E}^N}\setminus S \\
          I_j&  \textrm{ for } x\in I_j \textrm{ and } j\in J
                 \end{cases}
\ee
The {\bf smocked distance function}, $d: X\times X \to [0,\infty)$, is defined for $y, x\notin \pi(S)$, and intervals
$I_m$ and $I_k$ as follows:
\begin{eqnarray*}
d(\,x,\,y\,) &=& \min \left\{d_0(x,y), d_1(x,y), d_2(x,y), d_3(x,y), ...\right\} \\
d(\,x,\, I_k) &=& \min \{ d_0(x,z),  d_1(x,z), d_2(x,z), d_3(x,z), ...:\, z\in I_k\} \\
d(I_m,I_k) &=& \min \{ d_0(z',z), d_1(z',z), d_2(z',z), d_3(z',z), ... \,:z'\in I_m,\, z\in I_k \}
\end{eqnarray*}
where $d_k$ are the sums of lengths of segments that jump to and between $k$ intervals:
\begin{eqnarray*}
d_0(v,w) &=& |v-w|\\
d_1(v,w) &=&  \min\{ |v-z_1|+|z'_1-w|:\, z_1, z_1'\in I_{j_1}, \, j_1 \in J\}\\
d_2(v,w) &=& \min\{ |v-z_1|+|z'_1-z_2|+|z'_2-w|:\, z_i, z'_i\in I_{j_i}, \, j_1\neq j_2 \in J\}\\
d_k(v,w) &=& \min \{  |v-z_1|+\sum_{i=1}^k |z'_i-z_i|+|z'_k-w|:\, z_i, z'_i\in I_{j_i}, \, j_1\neq \cdots \neq j_k \in J\}.
\end{eqnarray*}
We define the {\bf smocking pseudometric} $\bar{d}: {\mathbb{E}^N}\times {\mathbb{E}^N} \to [0, \infty)$
to be
$$
\bar{d}(v,w)= d(\pi(v), \pi(w))=\min \{d_k(v',w'): \,\pi(v)=\pi(v'),\, \pi(w)=\pi(w'),\, k\in {\mathbb{N}}\}.
$$
We will say the smocked metric space is {\bf parametrized by points in the intervals}, if 
\be\label{param-by-points}
J \subset {\mathbb{E}}^N \textrm{ and } \forall j \in J \,\, j \in I_j.
\ee
\end{defn}

In \cite{Smocked} it is proven that the minima are achieved:

\begin{thm} \label{thm-smock-metric}
The smocked metric space is a well defined metric space and in fact
for any $v,w  \in \mathbb{E}^N$ 
\be \label{exists-N}
\exists N(v,w) \le \lceil{|v-w|/\lambda}\rceil \,\, s.t.\,\, d_{N(v,w)}(v,w) \le d_{k}(v,w) \quad \forall k \ge {\mathbb{N}}
\ee
so the minimum in the definition of the smocking distance and pseudometric is achieved 
\be
\bar{d}(v,w)=d_N(v,w) \textrm{ and } d(\pi(v), \pi(w))=d_N(v,w).
\ee
\end{thm}

In that paper the following constants are defined as well:

\begin{defn} \label{defn-smocking-depth}
The smocking depth, $h$, is defined to be
\be
h=\inf\{ r:\, \mathbb{E}^N \subset T_r(S)\} \in [0,\infty].
\ee
which by definition of tubular neighborhood, is 
\be
h=\inf\{ r:\, \forall x \in X \,\,\exists j \in J \, \exists z \in I_j \,\,s.t.\,\, |x-z|<r\}.
\ee
\end{defn}

\begin{lem} \label{lem-smocking-depth}
The smocking depth satisfies:
\be
h= \sup \{ D(x): \, x\in \mathbb{E}^N\}
\ee
where $D: {\mathbb{E}^N}\to [0,\infty)$ is the distance to the interval set:
\be
 D(x)= \min\{ |x-z|: \, z\in I_j, \, j \in J \}.
\ee
\end{lem}

\begin{lem}
The infimum in Definition~\ref{defn-smocking-depth} is
achieved as is the supremum in Lemma~\ref{lem-smocking-depth}.
\end{lem}

\begin{defn}\label{defn-smocking-lengths}
The smocking lengths are defined
\begin{eqnarray}
L_{min}&=& \inf \{ L(I_j): \, j \in J\} \in [0,\infty)\\
L_{max}&=& \sup \{ L(I_j): \, j \in J\} \in (0,\infty]
\end{eqnarray}
and if $L_{min}=L_{max}$ we call the the smocking length.
{\em if the $I_j$ are not intervals we can replace 
length with diameter}.
\end{defn}

\begin{defn} \label{defn-sep}
The smocking separation factor, $\delta=\delta_X$, is defined to be
\be
\delta_X = \delta_\diamond \le \min\left\{ |z -w|: \,\, z\in  I_j, \, w\in I_k, \, j\neq k \in J\right\}.
\ee
\end{defn}

\begin{lem}\label{lem-hplusL}
If a smocked metric space is parametrized by points in intervals as in (\ref{param-by-points}),
then
\be
{\mathbb{E}}^N \subset T_{h+L}(S)
\ee
where $S = \bigcup_{j\in J} I_j$ is the smocking set and $h$ is the smocking depth
and $L=L_{max}$ is the maximum smocking length.
\end{lem}

\subsection{Review of Balls in Smocked Metric Spaces}  

\begin{figure}[h]
\label{fig:smocked-balls}
\includegraphics[width=.2 \textwidth]{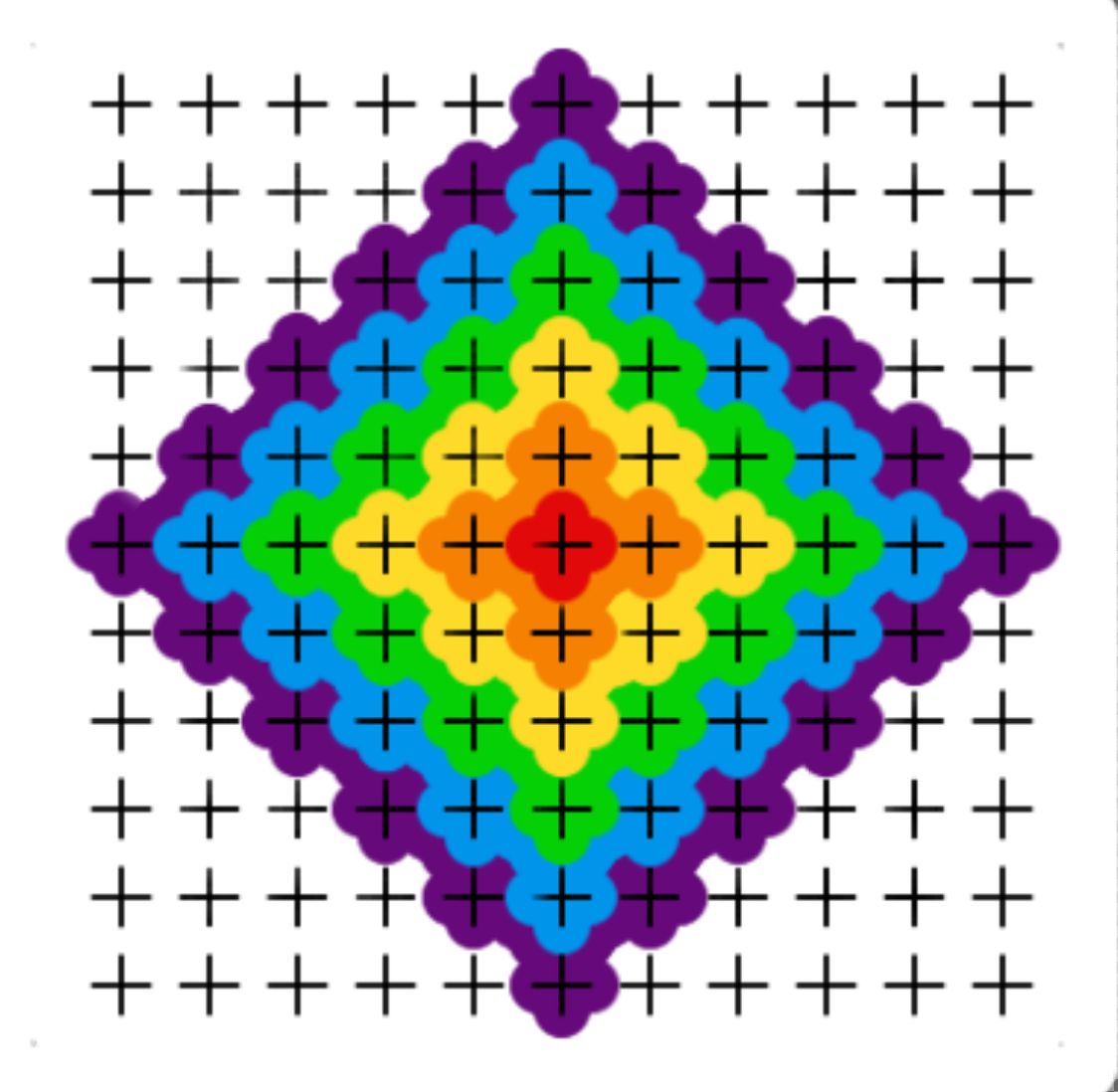} 
\includegraphics[width=.2 \textwidth]{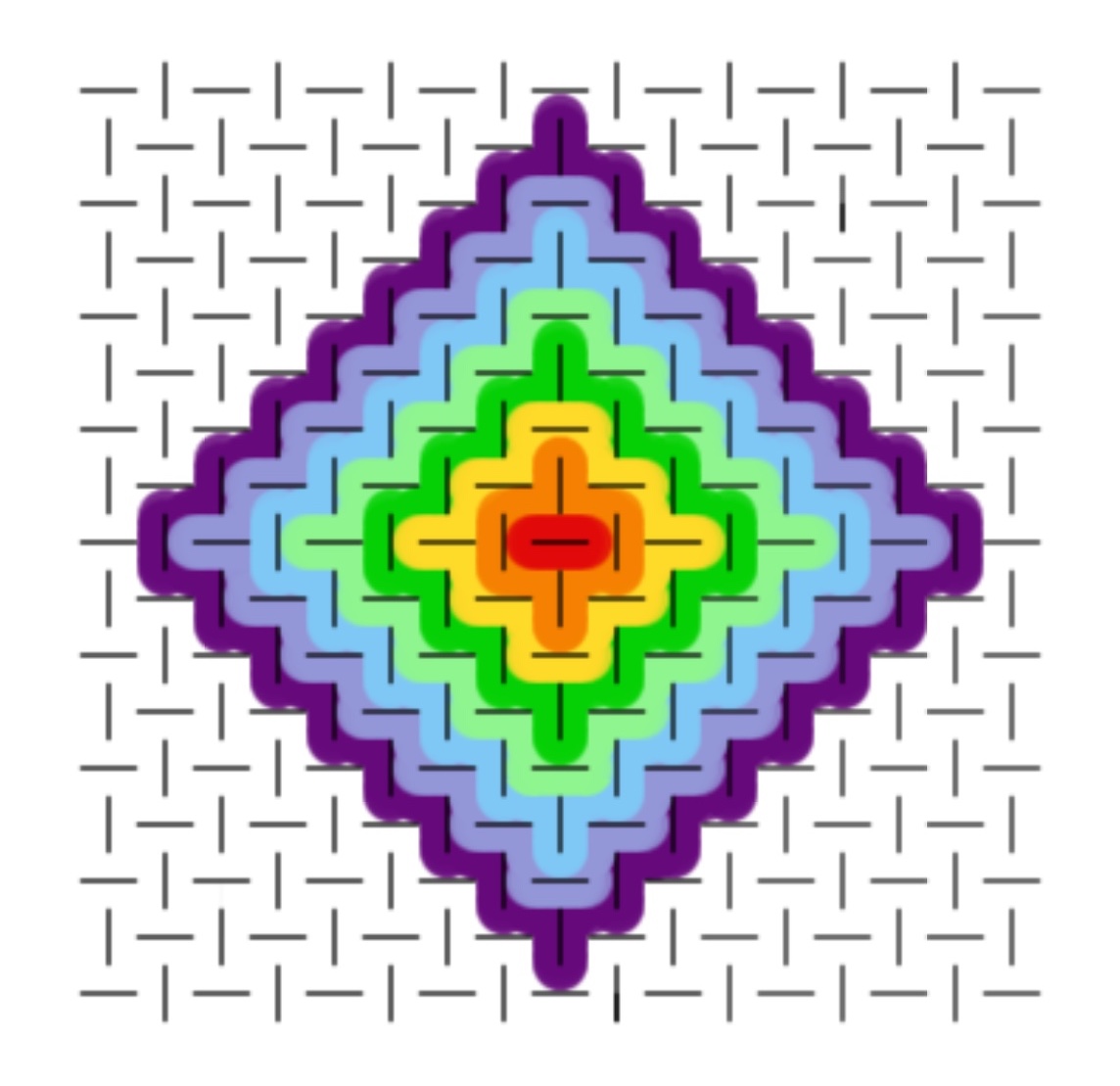} 
\includegraphics[width=.2 \textwidth]{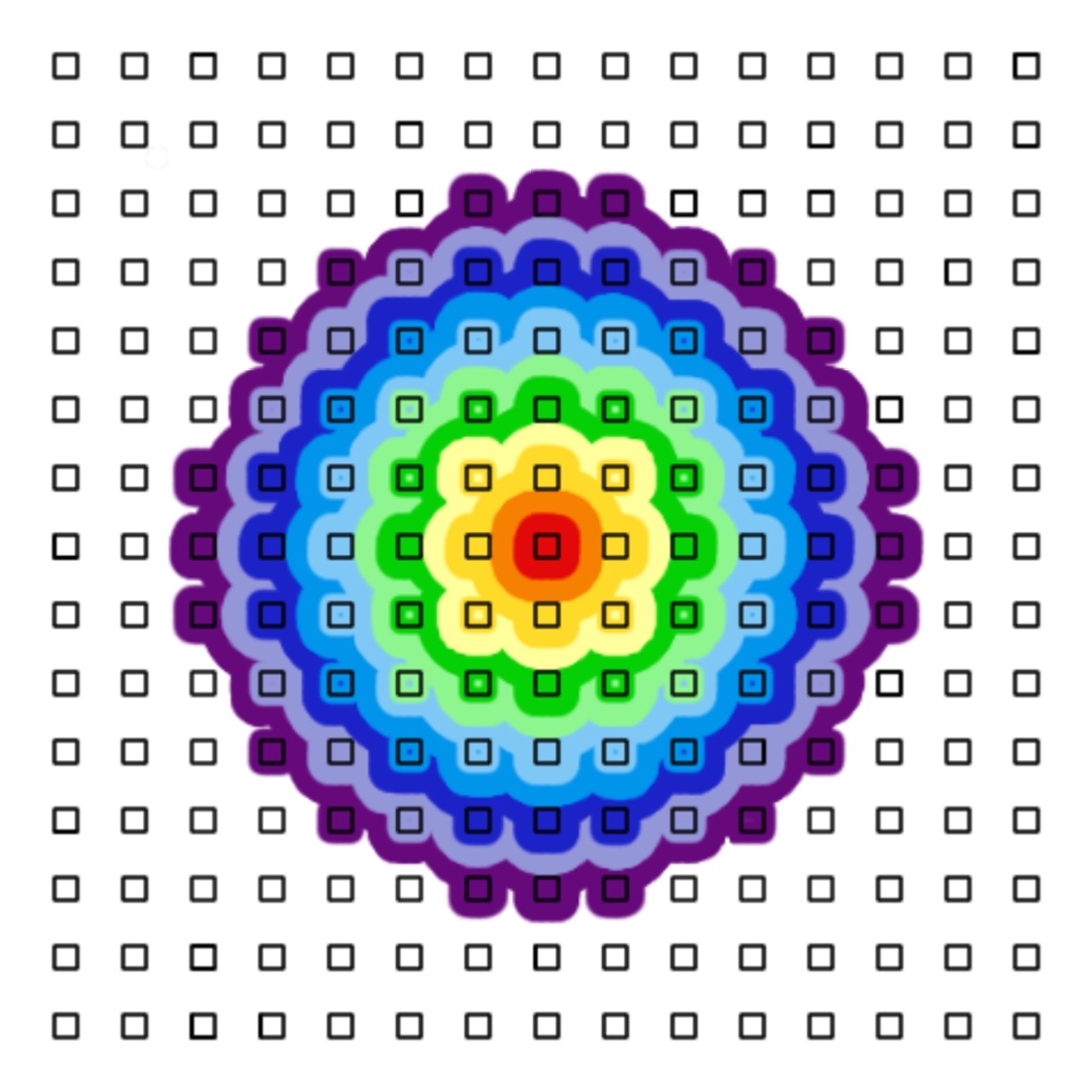} 
\caption{The balls in  $X_+$, $X_T$, and $X_\square$, 
were studied in \cite{Smocked}.}
\end{figure}

In \cite{Smocked}, the students explored the shapes of balls in a variety of smocked metric spaces.
To best describe these balls, one looks at their pre-images in Euclidean space: 
\be
U_r(x)=pi^{-1}(B_r(x)) \subset {\mathbb{E}}^N.
\ee
See Figure~\ref{fig:smocked-balls} for the balls found in that paper.  We will study the balls of
additional smocked spaces here using the following propositions and lemmas proven in that paper.

In \cite{Smocked} the following two propositions were proven:

\begin{prop} \label{ball-D(x)}
Suppose that $p\in {\mathbb{E}}^N\setminus S$ and $r<D(x)$ is defined as in Definition~\ref{defn-sep}, 
then
\be
\pi^{-1}(B_{r}(\pi(p)))=B_r(p)=\{x:\, |x-p|<r\}.
\ee
\end{prop}

\begin{prop} \label{ball-delta}
Suppose that $I_j$ is a smocking interval and $r<\delta_X$ is defined as in Definition~\ref{defn-sep}, 
then
\be
U_r(\pi(I_j))=\pi^{-1}(B_{r}(\pi(I_j)))=T_r(I_j)
\ee
\end{prop}

\subsection{Volumes and Smocked Metric Spaces}

In this paper we are finding the SWIF limit and SWIF limits
are defined using volumes and a more general notion called
mass which can be bounded using weighted volumes.  
In preparation for this, we will say something here
about the volumes of the pre-images of balls
which lie in Euclidean space:

\begin{lem}\label{Lip-ball-finite}
The volume, 
\be
\vol(U_R(x)) =\vol(\pi^{-1}(\bar{B}_R(x))
\ee
is finite and well defined because 
\be
\rho_x: {\mathbb{E}}^N \to [0,\infty) \textrm{ such that } \rho_x(v) = d_X(x, \pi(v))
\ee
is a proper Lipschitz function with respect to the Euclidean norm. 
\end{lem}

\begin{proof}
We first prove $\rho_x$ is Lipschitz:
\begin{eqnarray}
\frac{|\rho(v)-\rho(w)|}{|v-w|}&=&\frac{|d_X(x, \pi(v))-d_X(x, \pi(w))|}{|v-w|}\\
&\le &\frac{|d_X(\pi(v), \pi(w))|}{|v-w|} \le 1.
\end{eqnarray}
It is proper because its level sets are compact.
Since 
\be
U_R(x)= \rho_x^{-1}([0,R]) \textrm{ where } \rho_x(v) = d_X(x, \pi(v))
\ee
the volume may be computed exactly as in vector calculus by
breaking up the region into subregions with boundaries and
integrating the function $1$.
\end{proof}

\begin{lem}
If $\vare$ is less than the separation factor then
\be
\vol(U_\vare(\pi(I)))=\vol(T_\vare(I)) = \omega_N \vare^N + \omega_{N-1} \vare^{N-1} L.
\ee
where $\omega_N=\vol^{E}(B_1(0))$ and $\omega_0=2$.
\end{lem}

\begin{proof}
This follows from the fact that the tubular neighborhood is the union of a cylindrical region around interval and two hemispheres at the tips.
\end{proof}

Since we would like the volume of $U_R(x)$ to
relate well with the volume of $B_R(x)$
we now define a nice smocking set:

\begin{defn}\label{smocking-nice}
We say that a smocking set $S$ is nice
if for any fixed compact set, $K$, we have
\be
\lim_{\vare\to 0} \vol\left(T_{\vare}(S)\cap K\right) =0.
\ee
A smocked metric space is nice if it has a nice smocking set.
\end{defn}

Combining this definition with the two previous propositions
we immediately have the following new lemma:

\begin{lem}\label{nice-lim}
In a smocked space, $X$, with a nice smocking set,
\be
\lim_{\vare \to 0} \vol(U_\vare(x)) =0 \qquad \forall x\in X.
\ee
\end{lem}

\begin{rmrk}
Note that not all smocked metric spaces have nice smocking sets.  If for example
a smocking interval 
\be
I=[0,1]\times[0,1]\subset {\mathbb{E}}^2
\ee
then
\be
\lim_{\vare \to 0} \vol(U_\vare(\pi(I)))=\lim_{\vare \to 0} T_\vare(I) = \vol(I)=1\neq 0.
\ee
\end{rmrk}

\begin{lem}\label{who-nice}
A smocked metric space with $L_{max}<\infty$
has a nice smocking set
iff for every interval in the smocking set
\be
\lim_{\vare \to 0} \vol(T_\vare(I_j)) =0 \qquad \forall j\in J.
\ee
\end{lem}

\begin{proof}
If one interval has nonzero limit, just set $K$ to
be that interval to see that the smocking set is not
nice. Conversely, assume all the intervals
have $0$ limits and take any $K$.
Since the intervals in a smocking set $\{I_j: \, j\in J\}$
have a separation factor, and $K$ is compact,
there are only finitely many
intervals in the set
\be
\{I_j:\, j\in J_K\}=\{ I_j: \, I_j \cap K \ne \emptyset \,\, j\in J\}.
\ee
So we have a finite sum
\be
\vol(T_{\vare}(S)\cap K) \le  \sum_{j\in J_K} \vol(T_\vare(I_j))
\ee
Taking the limit we have our proof of the converse.
\end{proof}

\begin{rmrk} \label{all-four-nice}
It is easy to see that
 all four of our smocked spaces, $X_\diamond$,  $X_+$, $X_\square$, and $X_T$
are nice because the smocking intervals are line segments 
in $X_\diamond$ and $X_T$, they are unions of two
line segments in  $X_+$, and they are unions of four line segments in $X_\square$,
and the tubular neighborhood of a line segment, $I$, in ${\mathbb{E}}^N$ of length
$L$ has volume
\be
\vol(T_\vare(I)) = \omega_N \vare^N + \omega_{N-1} \vare^{N-1} L.
\ee
\end{rmrk}

\subsection{Review of GH Convergence and Tangent Cones at Infinity}

Gromov-Hausdorff convergence was first defined by Edwards in \cite{Edwards} and 
rediscovered by Gromov in \cite{Gromov-metric}.  

\begin{defn} \label{defn-GH}
We say a sequence of compact metric spaces
\be
(X_j, d_j) \GHto (X_\infty, d_\infty)
\ee
iff
\be
d_{GH}((X_j,d_j), (X_\infty, d_\infty)) \to 0.
\ee
Where the Gromov-Hausdorff distance is defined
\be
d_{GH}(X_j, X_\infty) = \inf \{d^Z_H(\varphi_j(X_j), \varphi_\infty(X_\infty)): \,\, Z,\,\, \varphi_j: X_j \to Z\}
\ee
where the infimum is over all compact metric spaces, $Z$, and over all
distance preserving maps $\varphi_j: X_j\to Z$:
\be
d_Z(\varphi_j(a), \varphi_j(b))=d_j(a,b) \,\,\,\forall a,b \in X_j.
\ee
The Hausdorff distance is defined 
\be
d_H(A_1, A_2) = \inf\{r:\,\, A_1\subset T_r(A_2) \textrm{ and } A_2\subset T_r(A_1) \}.
\ee
\end{defn}

Gromov also defined pointed Gromov-Hausdorff convergence:

\begin{defn}\label{defn-ptGH}
If one has a sequence of complete noncompact metric spaces, $(X_j, d_j)$, and points
$x_j \in X_j$, one can define pointed GH convergence:
\be
(X_j, d_j, x_j) \ptGHto (X_\infty, d_\infty, x_\infty)
\ee
iff
for every radius $R>0$ the closed balls of radius $R$ in $X_j$ converge in the GH sense as metric spaces
with
the restricted distance to closed balls in $X_\infty$:
\be
d_{GH}((\bar{B}_r(x_j)\subset X_j,d_j), (\bar{B}_r(x_\infty)\subset X_\infty, d_\infty)) \to 0.
\ee
\end{defn}

One may consider a single unbounded metric space, and take a sequence of rescalings of that
metric space.  A Gromov-Hausdorff limit of a sequence of rescalings, if it exists, is called a
GH tangent cone at infinity:

\begin{defn} \label{defn-tan-cone}
A complete noncompact metric space with infinite diameter, $(X, d_X)$,
has a GH tangent cone at infinity, $(Y, d_Y)$, if there is a sequence of
rescalings, $R_j \to \infty$, and points, $x_0\in X$ and $y_0\in Y$, such that
\be
(X, d/R_j, x_0) \ptGHto (Y, d_Y, y_0)
\ee
\end{defn}

There are a variety of theorems in the literature concerning the existence and
uniqueness of such tangent spaces at infinity.  

\subsection{Review of Tangent Cones at Infinity for Smocked Spaces}

\begin{figure}[h]
\label{fig:tan-cones}
\includegraphics[width=.8 \textwidth]{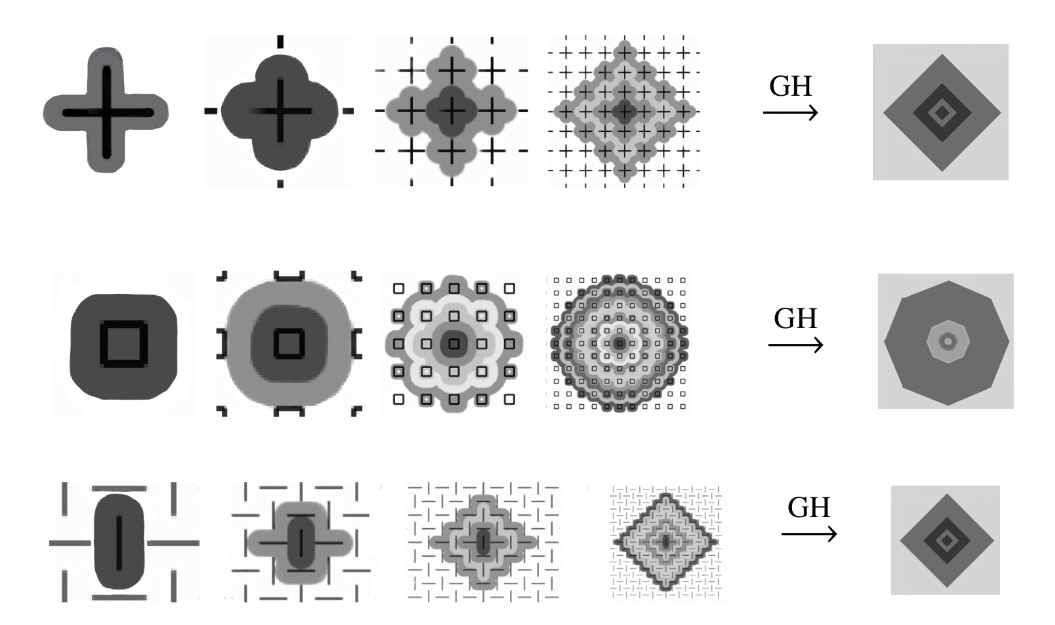}
%\includegraphics[width=.5 \textwidth]{diamondR} \raisebox{1 em}{$\quad{\GHto}\quad$}  
%\includegraphics[width=.13 \textwidth]{diamondL} \\  \quad \\
%\includegraphics[width=.5 \textwidth]{+R} \raisebox{1.9 em}{$\quad{\GHto}\quad$}   
%\raisebox{.5 em}{\includegraphics[width=.1 \textwidth]{taxiL} }\\  \quad \\
%\includegraphics[width=.5 \textwidth]{squareR}  \raisebox{1.8 em}{$\quad{\GHto}\quad$}  
%\raisebox{.5 em}{\includegraphics[width=.1 \textwidth]{octL}}\\ \quad \\
%\includegraphics[width=.5 \textwidth]{T-R} \raisebox{1.4 em}{$\quad{\GHto}\quad$}  
%\raisebox{.18 em}{\includegraphics[width=.08 \textwidth]{taxiL}}
\caption{In \cite{Smocked} the GH tangent cones at infinity were
found for  $X_+$, $X_\square$, and $X_T$ and proven to be normed spaces.}
\end{figure}

In \cite{Smocked} the GH tangent cones at infinity were found for four smocked
spaces.  See Figure~\ref{fig:tan-cones}.  It was shown that these four spaces had
unique tangent cones at infinity that were normed spaces.
The proofs of convergence were
based upon the lemmas and theorems stated in this subsection (which can
be applied more generally to a large class of smocked spaces to prove there
exist unique tangent cones at infinity which are normed spaces).  Note that
these results were proven using only the definitions by Gromov reviewed above.

The first step towards finding a tangent cone at infinity or any GH
limit of a smocked space is to apply the following lemma proven
in \cite{Smocked} to obtain an estimate on the smocking pseudometric
using only estimates on the distances between intervals.

\begin{lem} \label{lem-dil-to-approx}
Given an $N$ dimensional smocked metric space 
parametrized by points in intervals as in (\ref{param-by-points}),
with smocking depth, $h\in (0,\infty)$,
and smocking length $L=L_{max}\in (0,\infty)$, 
if one can find a Lipschitz
function $F: {\mathbb{E}}^N \to [0, \infty)$ such that
\be
|\,d(I_j, I_{j'})\,- \,[F(j)-F(j')] \,| \,\le \,C, 
\ee
then
\be
|\,\bar{d}(x, x')\,- \,[F(x)-F(x')] \,| \,\le \, 2h+ C + 2 \dil(F) (h+L)
\ee
where $dil(F)$ is the dilation factor or Lipschitz constant of $F$:
\be
\dil(F) = \sup \left\{\frac{|F(a)-F(b)|}{|a-b|}\,:\,\, a\neq b \in {\mathbb E}^N\right\}.
\ee
\end{lem}

This next theorem, also proven in \cite{Smocked}, 
is crucial for finding the GH limits of rescalings:

\begin{thm}\label{thm-smocking-R}
Suppose we have an $N$ dimensional smocked metric space, $(X, d)$, as in Definition~\ref{defn-smock}
 such that
\be
|\,\bar{d}(x, x')\,- \,[F(x)-F(x')] \,| \,\le \, K \qquad \forall x,x' \in {\mathbb{E}}^N
\ee
where $F: {\mathbb{E}}^N \to [0,\infty)$ is a norm.
Then $(X,d)$ has a unique GH tangent cone at infinity, $({\mathbb{R}}^N, d_F)$,
where
\be
d_F(x,x')=||x-x'||_F=F(x-x').
\ee
\end{thm}

\section{Balls in Smocked Metric Spaces are Integral Current Spaces}\label{CurrentSpaces}

Recall that any compact metric space, (X,d), endowed
with a countable collection of biLipschitz charts such that the
total volume of the images of all the charts is finite is an integer rectifiable current space
\cite{SW-JDG}.  If
the boundary of this space is also integer rectifiable, then it is an integral
current space.  Before we define this in more detail, we will look at
bi-Lipschitz maps onto closed balls in our smocked space.  We will then
find the right collection of biLipschitz charts for a closed ball in a smocked space, and finally
prove they are integral current spaces.  Before each theorem we will
provide the precise definition of what we are proving.

Throughout this section, as above,
\be
U_R(p)=\pi^{-1}(\bar{B}_R(p)) \subset {\mathbb E}^N,
\ee
is the pre-image of a closed ball under the smocking map.

\subsection{Finding bi-Lischitz maps}

\begin{defn}
A map $f: Y\to Z$ is bi-Lipschitz if it is a bijection with uniform 
upper bounds on $\dil(f)$ and $\dil(f^{-1})$.  That is, there exists
$\lambda>0$ such that
\be
\frac{1}{\lambda} \le \frac{d_Z(f(p), f(q))}{d_Y(p,q)} \le \lambda \qquad \forall p,q\in Y.
\ee
\end{defn}

\begin{rmrk}
Observe immediately that even though the smocking map $\pi: {\mathbb{E}}^N \to X$ 
is Lipschitz, it is not bi-Lipschitz
onto its image.  It is not even bijective.
Even if we remove the smocking set and study the restriction 
\be
\pi: {\mathbb{E}}^N\setminus S \,\,\to\,\, X \setminus \pi(S),
\ee
we only obtain a bijective map but it is still not bi-Lipschitz map.   One can see
this by taking $p_i,q_i$ arbitrarily close to the same interval in the smocking set,
but keeping $\pi^{-1}(p_i)$ and $\pi^{-1}(q_i)$ a definite distance $L/2$ apart:
\be
\frac{|\pi^{-1}(p_i)-\pi^{-1}(q_i)|}{d_X(p_i,q_i)} \to \infty.
\ee 
\end{rmrk}

\begin{lem} \label{lem-Lip}
If we avoid a tubular neighborhood of the smocking set,
\be
\pi: {\mathbb{E}}^N\setminus T_r(S) \to X \setminus \pi(T_r(S)),
\ee
then we do have a bi-Lipschitz map.   This map is still bi-Lipschitz
when further restricted to a compact ball in $X$:
\be
\pi: U_R(p) \setminus T_r(S) \to \bar{B}_R(p) \setminus \pi(T_r(S)).
\ee
In fact on $U_R(p)$ we have
\be
\min\left\{\frac{r}{\diam(U_R(p))}, 1\right\} \le \frac{d(\pi(p), \pi(q))}{|p-q|} \le 1 \qquad \forall p,q\in U_R(p) \setminus T_r(S)
\ee
\end{lem}

\begin{proof}
We already know that $\dil(\pi)= 1$.  
So we have a bi-Lipschitz map, unless
there are 
\be
p_i, q_i \in X \setminus \pi(T_r(S))
\ee 
such that
\be
\lim_{i\to \infty} 
\frac{|\pi^{-1}(p_i)-\pi^{-1}(q_i)|}{d_X(p_i,q_i)} =\infty.
\ee 
This can only happen if 
\be
\lim_{i\to \infty} d_X(p_i,q_i) = 0.
\ee
Thus there is an $N$ sufficiently large that
\be
 d_X(p_i,q_i) < r \qquad \forall i \ge N.
 \ee
 Since
 \be
 |p_i-z| >r \textrm{ and } |q_i-z'|>r \qquad \forall z,z'\in S,
 \ee
 we see the minimum in the Definition~\ref{defn-smock} is achieved by
 a direct segment from $\pi^{-1}(p_i)$ to $\pi^{-1}(q_i)$:
 \be
 d_X(\pi^{-1}(p_i), \pi^{-1}(q_i))= |(\pi^{-1}(p_i)-\pi^{-1}(q_i)| \qquad \forall i \ge N.
 \ee
 Thus 
\be
\lim_{i\to \infty} 
\frac{|\pi^{-1}(p_i)-\pi^{-1}(q_i)|}{d_X(p_i,q_i)} = 1
\ee 
which is a contradiction.  

In fact we know for any $p,q\in U_R(p) \setminus T_r(S)$ with $d(p,q)<r$ we have:
\be
\frac{|\pi^{-1}(p)-\pi^{-1}(q)|}{d_X(p,q)} =1.
\ee
If $p,q\in U_R(p) \setminus T_r(S)$ has $d(p,q)\ge r$
then we have
\be
\frac{|\pi^{-1}(p)-\pi^{-1}(q)|}{d_X(p,q)} \le \frac{\diam(U_R(p))}{r}.
\ee
\end{proof}

\subsection{Review of Integer Rectifiable Currents}

The following definitions are in Ambrosio-Kirchheim \cite{AK}.

\begin{defn}
We say $(f_0, f_1,...,f_N)$ is an $N$-tuple on a complete metric
space $X$ if $f_0$ is bounded and if all $f_i: X \to {\mathbb{R}}$
are Lipschitz.   Note that if 
\be
\varphi: A \subset {\mathbb{E}}^N \to X
\ee
is Lipschitz, then 
\be
f_j \circ \varphi: A \subset {\mathbb{E}}^N \to {\mathbb{R}}
\ee
 is also Lipschitz and is thus differentiable almost everywhere.
 \end{defn}

\begin{defn}
Given a precompact Borel set $A\subset {\mathbb{E}}^N $ and a 
Lipschitz map 
\be
\varphi: A \subset {\mathbb{E}}^N \to \varphi(A)\subset X
\ee
we can define $\varphi_*\lbrack A\rbrack$ , denoted $\varphi_*\lbrack A\rbrack$,
which acts on an $N$-tuple:
\be
\varphi_*\lbrack A\rbrack(f_0, f_1,...,f_N)=\int_A f_0\circ \varphi \,\,d(f_1 \circ \varphi)\wedge \cdots \wedge
d(f_N \circ \varphi)
\ee
This integration is well defined because Lipschitz functions on Euclidean space
are differentiable almost everywhere.  
\end{defn}

\begin{defn}\label{defn-int-rect}
Given a countable collection of 
weights $a_i \in {\mathbb{Z}}$ and bi-Lipschitz charts
from Borel sets $A_i$ in Euclidean space of
dimension $N$, 
\be\label{charts}
\varphi_i: A_i \subset {\mathbb{E}}^N \to \varphi_i(A_i) \subset X
\ee
that are pairwise disjoint
\be\label{pairwise-disjoint}
\varphi_i(A_i)\cap \varphi_j(A_j) =\emptyset
\ee
whose  total weighted volume is finite
\be\label{weighted1}
\sum_{i=1}^\infty |a_i| \Lip(\varphi_i)^N \vol(A_i) < \infty,
\ee
we may define an $N$ dimensional {\bf integer rectifiable current}, $T$,
acting on $N$-tuples
\be\label{int-rect}
T(f_0, f_1,...,f_N)=\sum_{i=1}^\infty \varphi_{i*}\lbrack A_i\rbrack(f_0, f_1,...,f_N).
\ee
Two collections of weighted charts define the same
current if they have the same values when acting on $N$ tuples.
\end{defn}

\subsection{Review of the Mass of Currents}

Ambrosio-Kirchheim define a mass for an integer rectifiable current
which they then prove to be bounded above by a constant multiple of the 
weighted volume in Section 9 of \cite{AK}:
\be \label{mass}
\mass(T) \le C_N \sum_{i=1}^\infty |a_i| \mathcal{H}^N(\varphi_{i}(A_i)) 
\ee
where
\be \label{mass-constant}
C_N=2^N/ \omega_N.
\ee
The mass is finite of the currents we've defined here because
\be 
\mathcal{H}^N(\varphi_{i}(A_i)) \le \Lip(\varphi_i)^N \vol(A_i).
\ee
We do not need the precise definition of mass in this paper.

\begin{rmrk}\label{push}
Ambrosio-Kirchheim defined the push forward of a current
by a Lipschitz map $\psi: X \to Y$ to be
\be
\psi_\#( T) (f_0,f_1,...,f_N) = T(f_0\circ \psi, f_1 \circ \psi,...,f_N\circ \psi)
\ee
and prove that
\be \label{push-mass}
\mass(\psi_\# T) \le (\dil \psi)^N \mass(T).
\ee 
\end{rmrk}

\subsection{Review of Adding and Subtracting Currents}

Ambrosio-Kirchheim defined addition of currents,
\be
(T+S)(f_0, f_1,...,f_N)= T(f_0, f_1,...,f_N) + S(f_0, f_1,...,f_N).
\ee
One can prove that when $T,S$ are integer rectifiable then so is
their sum (by carefully cancelling parts of the various charts in the
union of all the charts as necessary to ensure they are pairwise
disjoint).  They prove that
\be
\mass(T+S) \le \mass(T) + \mass(S)
\ee
Note this is not an equality as some charts in the parametrization
of $T$ may cancel with some charts in the parametrization of $S$.
It is easy to verify that
\be
\partial(T+S)=\partial T + \partial S.
\ee
Thus the sums of integral currents is an integral current. 
The $0$ current has
\be
0(f_0,...,f_N)=0 \qquad \forall \textrm{ tuples } (f_0,...,f_N).
\ee
Note that the negative of a current can be parametrized 
by the same collection of charts as the original current
but with the opposite orientation on all of them:
\begin{eqnarray*}
-T(f_0, f_1,...,f_N)&=&\sum_{i=1}^\infty -\varphi_{i*}\lbrack A_i\rbrack(f_0, f_1,...,f_N)\\
&=& \sum_{i=1}^\infty -\int_{A_i} (f_0\circ\varphi_i)\, d(f_1\circ \varphi_i)\wedge \cdots \wedge (f_N\circ \varphi_i)\\
\end{eqnarray*}

\subsection{Review of Integer Rectifiable Current Spaces}

The following definition first appeared in \cite{SW-JDG} by Sormani-Wenger:

\begin{defn}
Suppose a complete metric space, $(X,d)$, is endowed
with a countable collection weights $a_i \in {\mathbb{Z}}$
(here $a_i=1$) and bi-Lipschitz charts from Borel sets, $A_i$, in Euclidean space of
dimension $N$, as in (\ref{charts})
that are pairwise disjoint as in (\ref{pairwise-disjoint})
such that total weighted volume is finite as in (\ref{weighted1})
and such that
 the complement of the sets has zero measure
\be\label{zero-measure}
\mathcal{H}^N\left( X \setminus \bigcup_{i=1}^\infty \varphi_i(A_i)\right) =0.
\ee
Then we say the space has an
{\bf integer rectifiable current structure},
 \be\label{T}
 T= \sum_{k=0}^\infty a_i \pi_* \lbrack A_k \rbrack,
 \ee
 with mass as in \ref{mass}.
 We say the space has weight 1 if all the $a_i=1$.
 We say that $(X', d, T)$ is an {\bf integer rectifiable current space} where $X'\subset \bar{X}$
 is defined as the set of positive density of $T$ (here we need only know $X'$
 is defined).  
\end{defn}

\subsection{Balls in Normed Spaces are Integer Rectifiable Current Spaces}

%MOSHE AND HINDY DID THIS SECTION

Consider a metric space which is a finite dimensional normed vector space,
$({\mathbb{E}}^N, ||\cdot||_F)$ and compare it to Euclidean space.
\begin{rmrk}\label{equiv_norms}
We now use the fact that all norms on a finite dimensional vector space are equivalent. In other words for every norm $F$ on $\mathbb{E}^n$ there exists $C$ so that 
\be
\frac{1}{C}\leq \frac{|v|}{||v||_F}\leq C
\ee
for all $v$ in $\mathbb{E}^n$
\end{rmrk}

\begin{lem}\label{lambda}
There exists
\be
\lambda= \max_{v\neq 0} \frac{|v|}{||v||_F} <\infty.
 \ee
 Thus
 \be
 \bar{B}^F_R(0)=\{v:\, ||v||_F \le  R\} \subset B_{\lambda R} (0)
 \ee
\end{lem}

\begin{proof}
For the first part, note that  
\be 
\lambda= \max_{v\neq 0} \frac{|v|}{||v||_F}\leq C <\infty
\ee
This gives us that $\frac{|v|}{||v||_F}\leq \lambda$, so $|v|\leq \lambda ||v||_F$.
Supposing $v$ has $||v||_F\leq R$, we get
\be
|v| \leq \lambda ||v||_F
\ee
\end{proof}

\begin{thm}\label{thm-normed-ball}
A closed ball $\bar{B}^F_R(0)$ in an oriented finite dimensional 
normed vector space with norm $||x||_F$
is an integral current space with weight one where
the integral current structure, $T$, is defined by a single chart that is the identity
map and the ball itself as the domain of the chart.  Furthermore
\be
\mass(T_R) \le (\dil(id))^N \omega_N R^N\le C^N \omega_N R^N
\ee
and
\be
\mass(\partial T_R) \le (\dil(id))^N \omega_N N R^{N-1}\le C^N \omega_N N R^{N-1}
\ee
\end{thm}

\begin{proof}

The identity map 
\be 
id: \, (\bar{B}^F_R(0), d_E) \to (\bar{B}^F_R(0), d_F)
\ee
where $d_E(v,w)=|v-w|$ and $d_F(v,w)=||v-w||_F$
is biLipschitz 
with
\be
\dil(id) \le C < \infty
\ee
and
\be
\dil(id)^{-1} \le C \leq \infty
\ee
(Note that $\dil(id)^{-1}=\max_{v\neq 0} \frac{|v|}{||v||_F}=\lambda$ and $\dil(id) = \max_{v\neq 0} \frac{||v||_F}{|v|}$. By remark \ref{equiv_norms} we know there exists some C which bounds both $\max_{v\neq 0} \frac{|v|}{||v||_F}$ and $\max_{v\neq 0} \frac{||v||_F}{|v|}$).

So we define
\be
T_R=id_\# \lbrack \bar{B}^F_R(0) \rbrack
\ee
and this has
\be
\mass(T_R) \le (\dil(id))^N \vol_E(\bar{B}^F_R(0)) < (\dil(id))^N \omega_N R^N
\ee
by Lemma~\ref{lambda}.

\end{proof}

\subsection{Balls in Smocked Spaces are Integer Rectifiable Current Spaces}

We now prove that a closed ball in a smocked metric space is an integer rectifiable space
providing a precise set of charts to define a canonical integer rectifiable current on the space
(up to sign/orientation):

\begin{thm} \label{T-B-R}
A closed ball, $\bar{B}_R(p)$,
 in a smocked metric space, $(X, d_X)$, with smocking map $\pi: {\mathbb{E}}^N \to X$, 
 has a pair of
natural integer rectifiable current structures of weight $1$ defined by pushing forward
the two oriented local integral current structures, $\pm \lbrack U_R(p) \rbrack$, on ${\mathbb{E}}^N$: 
\be \label{Ak}
T=  \pi_\# \lbrack U_R(p) \rbrack = \sum_{i=0}^\infty \pi_* \lbrack A_i \rbrack
\ee
where
\be
A_k = U_R(p) \cup {T}_{1/k}(S) \setminus T_{1/(k+1)}(S) \textrm{ and } 
A_0= U_R(p) \setminus \left( S \cup \bigcup_{k=1}^\infty A_k\right)
\ee
so that for any collection of functions we have
\be
T(f, h_1,...,h_N)= \sum_{k=0}^\infty \int_{A_k} (f\circ \pi)\, d(h_1\circ \pi)\wedge\cdots\wedge d(h_N\circ \pi).
\ee
Furthermore the mass satisfies $\mass(T)=\vol(U_R(p))$.
\end{thm}

\begin{proof}
By Lemma~\ref{lem-Lip} these charts are bi-Lipschitz.  They are pairwise disjoint since $\pi$ is
a bijection away from the smocking set and the $A_k$ are pairwise disjoint.
Observe that (\ref{Ak}) holds because 
\be
\mathcal{H}^N \left(  U_R(p)\setminus \bigcup_{k=0}^\infty A_k\right) = \mathcal{H}^N \left( S \cap U_R(P)\right)=0
\ee
because the smocking set $S$ has zero measure.  Note that this also implies that
we have (\ref{zero-measure}).  

We claim that
\be
\mathcal{H}^m(\pi(A_k))=\mathcal{H}^m(A_k).
\ee
Recall that the Hausdorff measure, $\mathcal{H}^m(A_k)$ is defined using small sets, $Z_w$ about $w\in A_k$.
Once $\diam(Z_w)<1/(4k)$, we have a isometries $\pi: Z_w\to \pi(Z_w)$.  The Hausdorff measure, $\mathcal{H}^m(\pi(A_k))$ is also defined using small sets, $Z_p$ about $p\in \pi(A_k)$.  Once $\diam(Z_p)<1/(4k)$, we have a isometries $\pi: Z_p\to \pi(Z_p)$.  Thus we are estimating the Hausdorff measures of these sets with isometric 
collections of small sets.

As a consequence, taking our weights, $a_k=1$, we have
weighted volume equal to mass:
\be\label{weighted2}
\mass(T)=\sum_{k=1}^\infty  \mathcal{H}^N(\pi(A_k))= \sum_{k=1}^\infty  \mathcal{H}^N(A_k)=\vol(U_R(p))<\infty.
\ee
\end{proof}

\subsection{Review of Boundaries of Currents}

Ambrosio-Kirchheim defined the boundary of a current as follows in \cite{AK}:

\begin{defn}\label{defn-bndry}
The boundary of an $N$ dimensional current to be the
following $(N-1)$ dimensional current:
\be 
\partial T(f_0, f_1,..., f_{N-1})= T(1, f_0,..., f_{N-1}).
\ee
\end{defn}

\begin{lem}\label{lem-stokes}
For nice enough sets, $A\subset {\mathbb{E}}^N$, (eg. with piecewise smooth boundary)
\be
\partial  \pi\#\lbrack A\rbrack = \lbrack \pi(\partial A) \rbrack
\ee
\end{lem}

\begin{proof}
By the definitions and Stoke's Theorem we have:
\begin{eqnarray*}
\partial  \pi_\#\lbrack A\rbrack (f_0, f_1,...,f_{N-1}) &=& \pi_\#\lbrack A\rbrack(1, f_0, f_1,...,f_{N-1}) \\
&=&\int_A 1 \, d(f_0\circ \pi)\wedge d(f_1\circ \pi) \wedge \cdots \wedge d(f_{N-1}\circ \pi) \\
&=&\int_A  d\left(\, ( f_0\circ \pi) \, d(f_1\circ \pi) \wedge \cdots \wedge d(f_{N-1}\circ \pi\, \right) \\
&=&\int_{\partial A} (f_0\circ \pi) \, d(f_1\circ \pi) \wedge \cdots \wedge d(f_{N-1}\circ \pi) \\
&=& \pi_\# \lbrack \partial A\rbrack (f_0, f_1,...,f_{N-1}).
\end{eqnarray*}
\end{proof}

\subsection{Boundaries of Balls in Smocked Spaces}

In general the boundary of a current can be difficult to compute since there is no reason for
the Borel sets to have nice boundaries.  However in the boundaries of the charts we found for
our smocking sets are very nice and cancel perfectly when the smocking set is nice as in 
Definition~\ref{smocking-nice}.

\begin{prop}\label{smock-bndry}
Suppose the smocking set, $S$, is nice as in Definition~\ref{smocking-nice}.
Then boundary of the current defined in Theorem~\ref{T-B-R} is
\be
\partial T= \pi_\# \lbrack \partial U_R(p) \rbrack
\ee
so that for any collection of functions we have
\be\label{eq-smock-bndry}
\partial T(f_0, f_1,...,f_N)=  \int_{\partial U_R(p)} (f_0\circ \pi)\, d(f_1\circ \pi)\wedge\cdots\wedge d(f_{N-1}\circ \pi).
\ee
Since $\partial U_R(p)$ is covered by a finite collection of Lipschitz maps $\varphi_i: H_i \subset {\mathbb{E}}^{N-1}
\to \partial U_R(p)$ which are used to define what we mean by the integral in (\ref{eq-smock-bndry}):
\be\label{eq-smock-bndry2}
\partial T(f_0, f_1,...,f_N)=  \sum_i \int_{H_i} (f_0\circ \pi\circ \varphi_i) \, d(f_1\circ \pi\circ \varphi_i)\wedge\cdots\wedge d(f_{N-1}\circ \pi\circ\varphi_i).
\ee
Thus $\partial T$ is an integer rectifiable current with charts
$\pi \circ \varphi_i: H_i \subset {\mathbb{E}}^{N-1}\to \partial B_R(p)$.
Furthermore the mass satisfies 
\be
\mass(\partial T)=\vol(\partial U_R(p)).
\ee   
\end{prop}

\begin{proof}
\begin{eqnarray*} 
\partial T(f_0, f_1,...,f_{N-1})&=&T(1, f_0, f_1,...,f_{N-1}) \\
&=& \sum_{k=0}^\infty \pi_* \lbrack A_i \rbrack (1, f_0, f_1,...,f_{N-1}) \\
&=& \sum_{k=0}^\infty \int_{A_i} 1 \, d(f_0\circ \pi)\wedge d(f_1\circ \pi) \wedge \cdots \wedge d(f_{N-1}\circ \pi) \\
&=&\sum_{k=0}^\infty\int_{A_k}  d\left(\, ( f_0\circ \pi) \, d(f_1\circ \pi) \wedge \cdots \wedge d(f_{N-1}\circ \pi\, \right) \\
&=&\sum_{k=0}^\infty\int_{\partial A_k} (f_0\circ \pi) \, d(f_1\circ \pi) \wedge \cdots \wedge d(f_{N-1}\circ \pi) \\
&=& \int_{\partial A_0} (f_0\circ \pi) \, d(f_1\circ \pi) \wedge \cdots \wedge d(f_{N-1}\circ \pi) \\
&&+ \sum_{k=1}^\infty  
\,\, \int_{\partial A_k\cap U_R(p)} (f_0\circ \pi) \, d(f_1\circ \pi) \wedge \cdots \wedge d(f_{N-1}\circ \pi) \\
&& \qquad +  \int_{\partial T_{1/k}(S)\cap U_R(p)} (f_0\circ \pi) \, d(f_1\circ \pi) \wedge \cdots \wedge d(f_{N-1}\circ \pi) \\
&& \qquad -  \int_{\partial T_{1/(k+1)}(S)\cap U_R(p)} (f_0\circ \pi) \, d(f_1\circ \pi) \wedge \cdots \wedge d(f_{N-1}\circ \pi) 
\end{eqnarray*}
Telescoping the second with the third parts of the sum and totaling, we have
\begin{eqnarray*} 
\partial T(f_0, f_1,...,f_{N-1})&=&
 \int_{\partial U_R(p)} (f_0\circ \pi) \, d(f_1\circ \pi) \wedge \cdots \wedge d(f_{N-1}\circ \pi) \\
&&\quad - \, \lim_{k\to \infty}  \int_{\partial T_{1/(k+1)}(S)\cap U_R(p)} (f_0\circ \pi) \, d(f_1\circ \pi) \wedge \cdots \wedge d(f_{N-1}\circ \pi). 
\end{eqnarray*}
The limit is zero because the integral in the limit is
\begin{eqnarray*}
&\le & \int_{T_{k+1}(S)\cap U_R(p)} d(f_0\circ \pi) \wedge d(f_1\circ \pi) \wedge \cdots \wedge d(f_{N-1}\circ \pi)\\
&\le& \prod_{i=0}^{N-1} \Lip(f_i\circ \pi) \, \vol\left(T_{1/(k+1)}(S)\cap U_R(p)\right) \to 0
\end {eqnarray*}
because $S$ has dimension $N-1$ and $S \cap U_R(p)$ has finitely many components of finite 
$N-1$ volume.   
\end{proof}

\subsection{Review of Integral Currents and Spaces}

Ambrosio-Kirchheim defined  an integral current as follows in \cite{AK}:

\begin{defn}\label{defn-AK-integral}
An $N$ dimensional integral current is an 
$N$ dimensional integer rectifiable current whose
boundary is integer rectifiable.  The $0$
current
\be
0(f_0, f_1,...,f_N)=0
\ee
is also included among the integral currents.
\end{defn}

This allowed Sormani-Wenger to naturally define an integral current space in \cite{SW-JDG}.

\begin{defn}{int-cur-space}
An $N$ dimensional {\bf integral current space} is an integer rectifiable metric space
whose integral current structure is an $N$ dimensional integral current.  We also
include $0=(\emptyset, 0,0)$ as an
integral current space of dimension $N$.
\end{defn}

\subsection{Balls in Normed Spaces are Integral Current Spaces}

Consider a metric space which is a finite dimensional normed vector space,
$({\mathbb{E}}^N, ||\cdot||_F)$ and compare it to Euclidean space.

\begin{thm}\label{thm-normed-ball-integral}
A closed ball $\bar{B}^F_R(0)$ in an oriented finite dimensional 
normed vector space with norm $||x||_F$
is an integral current space with weight one where
the integral current structure $T^F_R$ is defined by a single chart that is the identity
map and the ball itself as the domain of the chart.
Furthermore
\be
\mass (T^F_R) \le \kappa_F^{2N} \omega_N R^N
\ee
and
\be
\mass (\partial T^F_R) \le \kappa_F^{N-1} \alpha R^{N-1}
 \ee
 where $\kappa_F >0$ is defined so that
 \be
 \frac{1}{ \kappa_F} \le \frac{||v||_F}{|v|} \le \kappa_F \qquad \forall v\neq 0.
 \ee
 and $\alpha= \vol_E(F^{-1}(1))$.
\end{thm}

\begin{proof}
Recall that in Theorem~\ref{thm-normed-ball}
we defined its integral current structure
\be
T=id_\# \lbrack \bar{B}^F_R(0) \rbrack.
\ee
Here $id$ is the biLipschitz identitiy map 
\be
id: \, (\bar{B}^F_R(0), d_E) \to (\bar{B}^F_R(0), d_F)
\ee
where $d_E(v,w)=|v-w|$ and $d_F(v,w)=||v-w||_F$
is biLipschitz 
with
\be
\dil(id) = \sup_{v\neq w} \frac{ ||v-w||_F }{|v-w|}\le \kappa_F.
\ee
and
\be
\dil(id)^{-1} =\sup_{v\neq w} \frac{|v-w|}{ ||v-w||_F } \le \kappa_F.
\ee
Thus
\be
\mass(T) \le (\dil(id))^{N} \vol_E(F^{-1}[0,R]) \le \kappa_F^{N} \omega_N(\kappa_F R)^N
\ee
because
\be
F^{-1}[0,R]=\{v\,|\, ||v||_F \le R\} \subset \{v\,| \,  |v| \le \kappa_F R\}.
\ee
Finally
\be
\mass(\partial T) \le (\dil(id))^{N-1} \vol_E(F^{-1}(\{R\}) \le \kappa_F^{N-1} \alpha R^{N-1}
\ee
because the volumes of level sets of a norm scale as follows:
\be
\vol_E(F^{-1}(\{R\})=\vol_E(F^{-1}(\{1\}) R^{N-1}.
\ee
\end{proof}

\subsection{Balls in Smocked Metric Spaces are Integral Current Spaces}

\begin{thm}\label{thm-integral}
Suppose the smocking set, $S$, is nice as in \ref{smocking-nice}.
Then the integer rectifiable current space, $(\bar{B}_R(p), d_X, \pi_*(U_R(p)))$,
defined in Theorem~\ref{T-B-R} is an integral current space.   Furthermore
the rescaled ball $(\bar{B}_R(p), d_X/t, \pi_*(U_R(p)))$ is also an integral current
space.
\end{thm}

\begin{proof}
The first part follows from the definition and Proposition~\ref{smock-bndry}.
The second part follows from the observation that any collection of charts
which is bi-Lipschitz with respect to $d_X$ is also bi-Lipschitz with respect to
$d_X/t$ and that the weighted volume will still be bounded.  
\end{proof}

\begin{rmrk}
It is an open question as to whether this is true or false for arbitrary smocked metric spaces.
See the end of the proof of Proposition~\ref{smock-bndry} to see where the hypothesis
on the smocking set was applied.
\end{rmrk}

\section{SWIF Convergence of Smocked Metric Spaces} \label{SWIF}

In this section we prove Theorem~\ref{thm-SWIF=GH-R}.   
Before beginning we
quickly review the definition of SWIF convergence in one subsection.

\subsection{Review of SWIF Convergence}

The Sormani-Wenger Intrinsic Flat (SWIF) distance was defined in \cite{SW-JDG} imitating
the Gromov-Hausdorff (GH) distance replacing the Hausdorff distance in Gromov's infimum
with the Flat distance of Federer-Flemming.  Since the Federer-Flemming flat
distance was defined only for integral currents in Euclidean space, we used
the notion of an integral current defined as in Ambrosio-Kirchheim \cite{AK} that we have 
just reviewed above.

\begin{defn} \label{defn-SWIF}
We say a sequence of compact integral current spaces
\be
(X_j, d_j, T_j) \Fto (X_\infty, d_\infty, T_\infty)
\ee
iff
\be
d_{SWIF}((X_j,d_j, T_j), (X_\infty, d_\infty, T_\infty)) \to 0.
\ee
Where the Sormani-Wenger intrinsic flat distance is defined
\be
d_{SWIF}(X_j, X_\infty) = \inf \{d^Z_F(\varphi_{j\#}(T_j), \varphi_{\infty\#}(T_\infty)): \,\, Z,\,\, \varphi_j: X_j \to Z\}
\ee
where the infimum is over all complete metric spaces, $Z$, and over all
distance preserving maps $\varphi_j: X_j\to Z$:
\be
d_Z(\varphi_j(a), \varphi_j(b))=d_j(a,b) \,\,\,\forall a,b \in X_j.
\ee
The Flat distance between two integral currents in $Z$ is defined 
\be
d_F^Z(S_1, S_2) = \inf\{\mass(A) + \mass(B):\,\, S_1-S_2 = A +\partial B\, \}
\ee
where the infimum is over all integral currents $A, B$ in $Z$
such that
\be
S_1(f_0,f_1,...,f_N)-S_2(f_0,f_1,...,f_N)= A(f_0,f_1,...,f_N)+ B(1, f_0,f_1,...,f_N)
\ee
for all tuples $(f_0,f_1,...,f_N)$ on $Z$.
\end{defn}

\begin{rmrk}
Examples in \cite{SW-JDG} demonstrate that GH and SWIF limits need not agree
and that SWIF limits may exist when there is no GH limit.
\end{rmrk}

In \cite{SW-JDG}, Sormani and Wenger proved that:

\begin{thm}
If $(X_j, d_j, T_j)$ converge to $(X_\infty, d_\infty, T_\infty)$
in the Lipschitz sense:
\be
\exists \textrm{ bi-Lip } F_j: X_j \to X_\infty \textrm{ with } \tfrac{1}{\lambda_j} < \dil(F_j) < \lambda_j
\textrm{ where } \lambda_j \to 1
\ee
and if 
\be
F_{j\#} T_j = T_\infty
\ee
(which holds if they are all weight $1$ and oriented in the same way),
then
\be
(X_j, d_j, T_j) \Fto (X_\infty, d_\infty, T_\infty).
\ee
\end{thm}

\begin{rmrk}
Our rescaled smocked metric spaces do not converge in the Lipschitz
sense to their tangent cones at infinity.  In fact there does not even exist
a bi-Lipschitz map between a smocked metric space with a nonempty
smocking set and a Euclidean space endowed with a definite norm.  We 
won't prove this claim in general
but will observe that the estimate obtained
in (\ref{lem-dil-to-approx}):
\be
|\,\bar{d}(x, x')\,- \,[F(x)-F(x')] \,| \,\le \, 2h+ C + 2 \dil(F) (h+L)
\ee
 are not able to control ratios of distances
between pairs of points $x, x'$ which lie on a common smocking 
interval. 
\end{rmrk}

\subsection{The set up}
%Team set: Sormani Kazaras Moshe

Take any $x_0\in X$.  By shifting the smocking set, $S$, we may assume 
that $\pi(0)=x_0$ where
$\pi: {\mathbb{E}}^N\to X$ is the pulled thread map.   
Let
\be
\lambda= \max_{v\neq 0} \frac{|v|}{|F(v)|} <\infty.
 \ee
 be as in Lemma~\ref{lambda}.  

We need to show that for all $r>0$
\be
(\bar{B}^X_{Rr}(x_0), d_X/R) \,\, \to \,\,(\bar{B}^F_r(0), d_F))
\,\,\textrm{ as }\,\, R\to \infty
\ee
where
\be
\bar{B}^F_r(0)=F^{-1}([0,r])=\{v:\,\, F(v)\le r\} \subset \{v:\,\, |v|\le \lambda r\} \subset B_{\lambda r}(0) \subset \mathbb{E}^N.
\ee
and
\be
\bar{B}^X_{Rr}(x_0)=\{x\in X: d_X(x,x_0) \le Rr\} \subset X
\ee
Since we are rescaling this ball in $X$, using the metric $d_X/R$, there is an isometry
\be
F_R: \,\,(\bar{B}^X_{Rr}(x_0), d_X/R) \to \,(\bar{B}^{R}_r(x_0), d_R)
\ee
whose image is a ball is the rescaled smocked metric space:
\be
\bar{B}^{R}_r(x_0)= \{x\in X_R: d_{R}(x,x_0) \le r\} \subset X_R
\ee
where $(X_R, d_{R})$ is the rescaled smocked metric space
defined with rescaled smocking intervals:
\be
{\mathcal{I}}_R=\{ I_{j}/R: \, j\in J\} \textrm{ where } I_{j}/R=\{z/R: \, z\in I_j\}
\ee
and a rescaled smocking map, 
\be
\pi_R: {\mathbb{E}}^{N} \to X_R 
\ee
Observe that
\be
\pi_R(v)=F_R( \pi(Rv))
\ee
See Figure~\ref{fig:tan-cones} to see how the smocking intervals rescale
in a variety of smocked metric spaces as $R\to \infty$.

We set
\be
U^R_r(x_0)=\pi_R^{-1}(\bar{B}^R_{r}(x_0)) \subset \mathbb{E}^N.
\ee
Observe that for $K/R <r$ we have
\begin{eqnarray}
U^R_r(x_0)&=&\{v: \,\,d_{R}(\pi_R(v), x_0)\le r\}\\
&=&\{v: \,\,d_{X}(\pi(Rv), x_0)\le rR\}\\
&=&\{v:\,\, \bar{d}_X(Rv, 0)\le rR\}\\
&\subset &\{v:\,\, ||Rv||_F<rR + K\}\\
&=&\{v:\,\, ||v||_F < r +(K/R)\}\\
&\subset& \{v:\,\, |v| < \lambda(r + (K/R)) \} \subset B_{2\lambda r}(0)\subset \mathbb{E}^N.
\end{eqnarray}
Note that the current structure, $T_R$, is
\be
T_R=\pi_{\#} \lbrack  U^R_r(x_0)\rbrack
\ee
as in Theorem~\ref{T-B-R}
and 
\be
T_F= id_{\#} \lbrack F^{-1}([0,r])\rbrack
\ee 
where $id$ is the identity map as in Theorem~\ref{thm-normed-ball-integral}

To prove the theorem, we will show that as $R\to\infty$
\be
d_{GH}\left((\bar{B}^{R}_{r}(x_0), d_{R}), (F^{-1}([0,r]), d_F) \right) \to 0
\ee
and
\be
d_{SWIF}\left((\bar{B}^{R}_{r}(x_0), d_{R}, T_R), (F^{-1}([0,r]), d_F, T_F) \right) \to 0.
\ee
\subsection{Constructing a metric space, $Z$}
% Team-d Ajmair, Maziar, Sormani

\begin{figure}[h]
\label{fig:SWIF-Z}
\includegraphics[width=.6\textwidth]{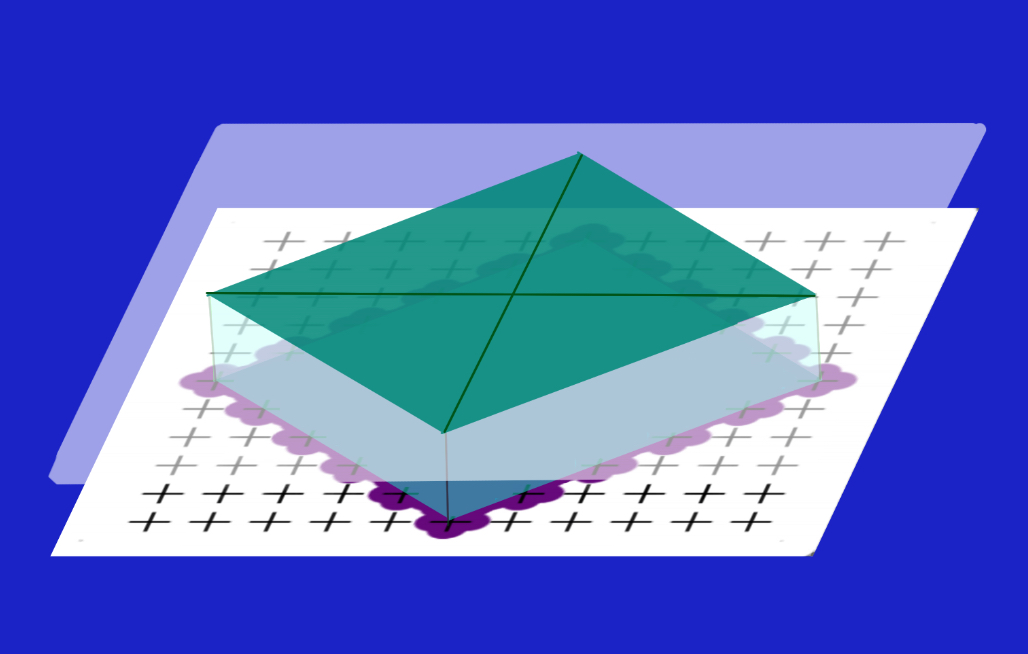} 
\caption{The metric space $Z$ is a smocked space 
defined using a collection of smocking intervals in the
plane ${\mathbb{E}}^N\times \{0\}$.  Also depicted here
are the images of our distance preserving maps.}
\end{figure}

As in
Figure~\ref{fig:SWIF-Z}, we take $Z$ to be a smocked space defined
in one dimension higher, ${\mathbb{E}}^{N+1}$, with
smocking intervals 
\be
{\mathcal{I}}_Z=\{ I_j/R\times\{0\}: \, j\in J\} \subset {\mathbb{E}}^{N} \times\{0\}
\ee
and a smocking map 
$
\pi_Z:\, {\mathbb{E}}^{N+1}\to Z.
$
and smocking metric 
\be
d_{Z}^S: Z \times Z \to [0, \infty).
\ee
Thus we have a distance preserving map
\be
\varphi^S: \, (X_R, d_R) \to (Z, d^S_Z)
\ee
defined by
\be \label{phiS}
\varphi^S (x) = \pi_Z ( \pi_R^{-1}(x) \times \{0\}).
\ee
 The smocking preimage of 
 \be
 \varphi^S(\bar{B}^R_r(x_0))=\pi_Z(U^R_{r}(0) \times \{0\})
 \ee
  is depicted in Figure~\ref{fig:SWIF-Z} as a 
 purple set in the lower plane with the intervals .
% ADD DETAILS?
 
\subsection{Choosing a height $H$:}
%Team-H Maziar Ajmain Sormani

We would like to show the map:
\be
\varphi^H: \, ({\mathbb{R}}^N , d_F ) \to \pi_Z( \, {\mathbb E}^N \times \{H\} )\subset Z
\ee
defined by
\be\label{phiH}
 \varphi^H(w) = \pi_Z(w,H)
 \ee
 to be distance preserving.  This will not work with $d_Z$, so we
 do not claim this is distance preserving.  For intuition, you
 will see the image of $F^{-1}([0,r])$ under this map as a green diamond in the upper plane in
 Figure~\ref{fig:SWIF-Z}.
 
 We choose $H=H_r>0$ such that
 \be\label{choice-H}
 H=\sqrt{ 8 \lambda r(K/R) + (K/R)^2}
 \ee
 so that  
 \be
 H\ge \sqrt{ 2a (K/R)+ (K/R)^2} \qquad \forall a\in [0,4\lambda r],
 \ee
 which implies
 \be \label{H-choice}
 \sqrt{H^2+a^2} \ge a+K/R \qquad \forall a\in [0,4 \lambda r].
 \ee
 Thus
 \be
 | (v,H) - (z,0) | \ge |v-z| +K/R \qquad \forall v,z \in B_{2\lambda r}(0) \subset {\mathbb E}^N
 \ee
 Since the smocking distance between points are sums
 over segments between intervals and all the intervals lie in ${\mathbb E}^N\times \{0\}$,
 for any $x, y \in B_{2\lambda r}(0)$ such that $|x-y|>\delta$ where $\delta$ is the smocking separation factor, there is $z \in B_{2\lambda}r(0)$ such that $(z,0)\in S_Z$ so that
 \begin{eqnarray}
 \bar{d}^S_Z( (x,H), (y,0)) &=& |(x,H)-(z,0)| +\bar{d}^S_Z((z,0),(y,0))\\
 &\ge &  |x-z| +K/R +\bar{d}^S_Z((z,0),(y,0))\\
 &\ge & \bar{d}^S_Z((x,0),(y,0)) + K/R\\
 &=& \bar{d}_X(x,y)+K/R \\
 &\ge & ||x-y||_F \textrm{ by the hypothesis rescaled}.
 \end{eqnarray} 
 Thus for all $x_i, y_i \in B_{2\lambda r}(0)\subset {\mathbb{E}}^N$
 \begin{eqnarray}
 \bar d_Z^S((y_1,H),(x_1,0)) + \bar d_Z^S((x_1,0),(x_2,0)) +\bar d_Z^S((x_2,0),(y_2,H))
& \geq&\\ ||y_1-x_1||_F + |x_1 - x_2| + ||x_2 - y_2||_F &\geq&\\ ||y_1-x_1||_F + ||x_1 - x_2||_F + ||x_2 - y_2||_F &\geq&\\ ||y_1- y_2||_F
 \end{eqnarray}
 and since all the smocking intervals lie in ${\mathbb E}^N\times \{0\}$
 and $||y_1-y_2||_F\le |y_1-y_2|$, we have
 \be  \label{Himplies}
 ||y_1-y_2||_F \le \bar{d}^S_Z( (y_1,H), (y_2,H)) \quad \forall y_1,y_2 \in B_{2\lambda r}(0)\subset {\mathbb{E}}^N.
 \ee
 Thus 
 \be 
\varphi^H: \, ({\mathbb{R}}^N , d_F ) \to  (Z, d_Z).
\ee
is distance nonincreasing.  Sadly it is not a distance
preserving map.

\subsection{Constructing a better metric space, $Z_R$.}
% TEAM-d Sormani Ajmain Maziar

We define a new metric space:
\be
Z_R=\pi_Z\left(B_{2\lambda r}(0) \times [0,H] \right)\subset Z
\ee
with a new metric
\be\label{dZR}
 d_Z^R(p_1,p_2)=\min\{ d_Z^S(p_1,p_2), d_Z^F(p_1,p_2)\}
 \ee
 where
 \be\label{d_Z^F def}
 d_Z^F(p_1,p_2)=
 \min\left\{ d_Z^S(p_1,\pi_Z(y_1,H))+ ||y_1-y_2||_F + d_Z^S(\pi_Z(y_2,H),p_2):\,\, y_i\in {\mathbb E}^N \right\}.
 \ee
Intuitively we are just shrinking the distances between points exactly enough
to make $\varphi^H$ a distance preserving function.  

\subsection{Proof that $d_Z^R$ is a metric on $Z_R$.}
% Team-d Ajmain and Maziar and Sormani

 It is easy to see that $d_Z^R$ is symmetric and definite because both
$d_Z^F$ and $d_Z^S$ are symmetric and 
\be
d_Z^S(x,y)=0 \,\,\iff\,\, x=y \,\,\textrm{ and }\,\,d_Z^F(x,y)=0\,\, \iff \,\,x=w=w'=y.  
\ee
 We claim the triangle inequality,
 \be
 d_Z^R(p_1,p_3)\le d_Z^R(p_1,p_2)+ d_Z^R(p_3,p_2) \qquad \forall p_1, p_2, p_3 \in Z_R.
 \ee
 {\em Case I:} Suppose both $d_Z^R(p_i,p_2)=d_Z^S(p_i,p_2)$ for $i=1,3$.
 Then the triangle inequality holds
 because $d_Z^S$ is a metric and 
 \be
 d_Z^R(p_1,p_2)\le d_Z^S(p_1, p_2).
 \ee
 {\em Case II:} Suppose only one is the smocking length:
 \be
 d_Z^R(p_1,p_2)=d_Z^F(p_1,p_2) \textrm{ and } d_Z^R(p_3,p_2)=d_Z^S(p_3,p_2)
 \ee
 Then there exists $y_i \in {\mathbb E}^N$
 such that
 \be
 d_Z^R(p_1,p_2)=d_Z^S(p_1,\pi_Z(y_1,H))+ ||y_1-y_2||_F + d_Z^S(\pi_Z(y_2,H),p_2).
 \ee
 Then we can apply these $y_i$  to estimate:
 \begin{eqnarray*}
 d_Z^R(p_1,p_3) &\le& d_Z^F(p_1,p_3)\\ 
 &\le& d_Z^S(p_1,\pi_Z(y_1,H))+ ||y_1-y_2||_F + d_Z^S(\pi_Z(y_2,H),p_3).\\
 &\le & d_Z^S(p_1,\pi_Z(y_1,H))+ ||y_1-y_2||_F + d_Z^S(\pi_Z(y_2,H),p_2)+d_Z^S(p_2,p_3)\\
 &=&d_Z^R(p_1,p_2) + d_Z^R(p_2,p_3).
 \end{eqnarray*}
 {\em Case III:} Suppose neither is the smocking length
 \be
 d_Z^R(p_1,p_2)=d_Z^F(p_1,p_2) \textrm{ and } d_Z^R(p_3,p_2)=d_Z^F(p_3,p_2)
 \ee
 Then there exists $y_{i,j} \in {\mathbb E}^N$
 such that
 \be
 d_Z^R(p_i,p_2)=d_Z^S(p_{i},\pi_Z(y_{i,i},H))+ ||y_{i,i}-y_{i,2}||_F + d_Z^S(\pi_Z(y_{i,2},H),p_2).
 \ee
 Then we can apply these $y_{1,1}, y_{3,3}$  to estimate:
 \begin{eqnarray*}
 d_Z^R(p_1,p_3) &\le& d_Z^F(p_1,p_3)\\ 
 &\le& d_Z^S(p_1,\pi_Z(y_{1,1},H))+ ||y_{1,1}-y_{3,3}||_F + d_Z^S(\pi_Z(y_{3,3},H),p_3).\\ 
 \end{eqnarray*}
 Since $||\cdot||_F$ is a norm, we have
 \be
 ||y_{1,1}-y_{3,3}||_F \le ||y_{1,1}-y_{1,2}||_F+||y_{1,2}-y_{3,2}||_F+||y_{3,2}-y_{3,3}||_F.
 \ee
 Combining this with the previous two equations we have
 \begin{eqnarray*}
 d_Z^R(p_1,p_3) &\le& d_Z^R(p_1,p_2) -d_Z^S(\pi_Z(y_{1,2},H),p_2) +||y_{1,2}-y_{3,2}||_F\\
&& + d_Z^R(p_3,p_2) -d_Z^S(\pi_Z(y_{3,2},H),p_2).
\end{eqnarray*}
So we need only show
\be
||y_{1,2}-y_{3,2}||_F\le d_Z^S(\pi_Z(y_{1,2},H),p_2)+d_Z^S(\pi_Z(y_{3,2},H),p_2) .
\ee
This holds by the norm being less than the
smocking length in (\ref{Himplies}) and the triangle inequality for $d_Z^S$.
Thus $(Z_R, d_Z^R)$ is a metric space.

\subsection{Proof that $\varphi_H$ of (\ref{phiH}) is now distance preserving.}
%Team-d Ajmain and Maziar

Let $\varphi_H$ be the map 
\be
\varphi_H: \, ( B_{2\lambda r}(0) , d_F ) \to (Z_R, d_{Z}^R)
\ee
defined by
\be\label{phi_H}
 \varphi_H(w) = \pi_Z(w,H)
 \ee
Note that the image of $\varphi_H$ is depicted as a green diamond shaped region
in the upper plane in Figure~\ref{fig:SWIF-Z}.

% Ajmain and Maziar Team-d
We want to show for any $x_1,x_2 \in B_{2\lambda r}(0) \subset \mathbb{E}^N$
\be
d_F(x_1,x_2)= d_Z^R(\varphi_H(x_1),\varphi_H(x_2)).
\ee
\begin{enumerate}
    \item["$\leq$"] 
        Equation \ref{Himplies} gives
        \begin{eqnarray}\label{phiH1}
            d_F(x_1,x_2) &\leq& d_Z^S(\varphi_H(x_1),\varphi_H(x_2)) .
        \end{eqnarray}
        All that remains to show is that 
        \begin{eqnarray}\label{phiH2}
            d_F(x_1,x_2) &\leq&d_Z^F(\varphi_H(x_1),\varphi_H(x_2)) .
        \end{eqnarray}
        By triangle inequality and equation \ref{phiH1} we have that for any $y_1, y_2 \in B_{2\lambda r}(0)$, 
        \begin{eqnarray}\label{almost phiH2}
            d_F(x_1,x_2) &\leq& d_F(x_1,y_1) + d_F(y_1,y_2) + d_F(y_2,x_2) \\
            &\leq& d_Z^S(\pi_Z(x_1,H), \pi_Z(y_1,H)) +  ||y_1 -y_2 ||_F \\
            && + \  d_Z^S( \pi_Z(y_2,H), \pi_Z(x_2,H)).
        \end{eqnarray}
        By passing to minimum of the right hand side of inequality \ref{almost phiH2} over all $y_i\in {\mathbb E}^N$, we get \ref{phiH2}, as required.
        
    \item["$\geq$"] By the definition of $d_Z^R$ we have
        \begin{eqnarray*}
            d_Z^R(\varphi_H(x_1),\varphi_H(x_2)) &=& d_Z^R(\pi_Z(x_1,H),\pi_Z(x_2,H))\\
            &=& \min\{ d_Z^S(\pi_Z(x_1,H),\pi_Z(x_2,H)), d_Z^F(\pi_Z(x_1,H),\pi_Z(x_2,H))\\
            &\leq& d_Z^F(\pi_Z(x_1,H),\pi_Z(x_2,H))\\
            &\leq& ||x_1 - x_2||_F = d_F(x_1,x_2)
        \end{eqnarray*}
    where the last step follows from taking $y_1=x_1$ and $y_2=x_2$ in definition \ref{dZR}.
\end{enumerate}

\subsection{Proof that $\varphi_S$ of (\ref{phiS}) is distance preserving.}
%Team-d Ajmain and Maziar

Let $\varphi_S$ be the map 
\be
\varphi_S: \, (\bar{B}^R_r(x_0) , d_R ) \to (Z_R, d_{Z}^R)
\ee
defined by
\be\label{phi_S}
 \varphi_S(w) = \pi_Z(\pi_R^{-1}(w)\times\{0\})
 \ee

Note that the image of $\varphi_S$ is depicted as a purple region
in the lower plane in Figure~\ref{fig:SWIF-Z}. Recall that
$\varphi^S$ was distance preserving with respect to $d_Z^S$.

Therefore, by definition \ref{dZR}, it suffices to show that \be
d_Z^S(\varphi_S(w_1),\varphi_S(w_2))\leq d_Z^F(\varphi_S(w_1),\varphi_S(w_2)).
\ee

\begin{figure}[h!]
\centering
\includegraphics[scale=.35]{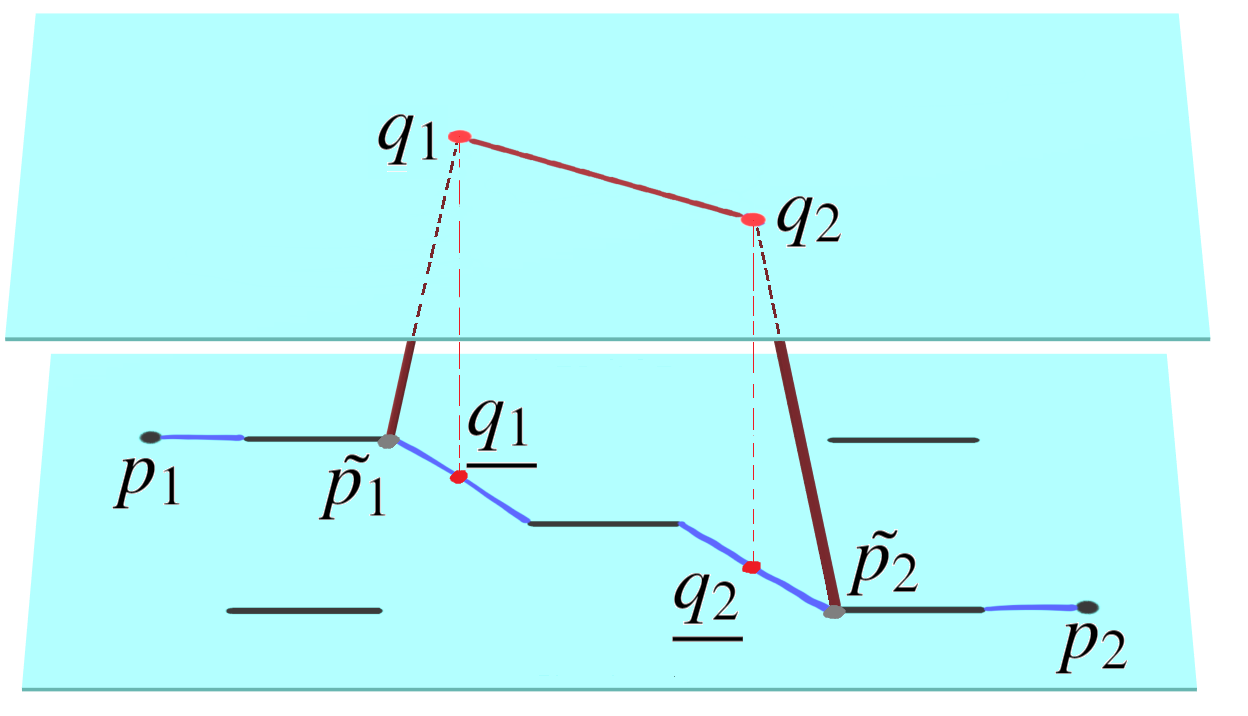}
\caption{The figure for Proof that $\varphi_S$ of (\ref{phiS}) is distance preserving.}
\label{dZF.png}
\end{figure}

Let $w_i = \pi_R(x_i)\in \bar{B}^R_r(x_0)$ for some $x_i\in \mathbb{E}^N$.  Let $p_i = \pi_Z(x_i, 0)$.
We wish to show 
\be
d_Z^S(p_1,p_2)\leq d_Z^F(p_1,p_2).
\ee
Intuitively, it means that to go from $p_1$ to $p_2$, it is better to stay on the $0$-plane rather than jumping to the $H$-plane and coming back to the $0$-plane (Fig. \ref{dZF.png}).

By definition of $d_Z^F$ in equation \ref{d_Z^F def}, there exists $y_1, y_2 \in \mathbb{E}^N$ such that 
\be \label{dZ_F(p_1,p_2)}
d_Z^F(p_1,p_2) =  d_Z^S(p_1,q_1)+ ||y_1-y_2||_F + d_Z^S(q_2,p_2)
\ee
where $q_i= \pi_Z(y_i,H)$.  Write $\underline{q_i}=\pi_Z(y_i,0)$.

Let $\tilde{p_i}=\pi_Z(\tilde{x_i},0)$ be points such that 
\be \label{d_Z^S(p_i, q_i)}
d_Z^S(p_i, q_i)= d_Z^S(p_i, \tilde{p_i}) + |(\tilde{x_i},0)-(y_i,H)|.\ee

By triangle inequality, 
    \be d_Z^S(p_1,p_2)\leq d_Z^S(p_1,\tilde{p_1})+d_Z^S(\tilde{p_1},\tilde{p_2})+d_Z^S(\tilde{p_2},p_2).\ee

By equation \ref{dZ_F(p_1,p_2)} and \ref{d_Z^S(p_i, q_i)}
\begin{eqnarray}
    d_Z^F(p_1,p_2)&=& d_Z^S(p_1,q_1)+||y_1-y_2||_F+d_Z^S(q_2,p_2)\\
    &=&d_Z^S(p_1,\tilde{p_1})+ |(\tilde{x_1},0)-(y_1,H)|+
    ||y_1-y_2||_F+  \\& & |(y_2,H)-(\tilde{x_2},0)|+ d_Z^S(\tilde{p_2},p_2) \nonumber
\end{eqnarray}

Therefore, it is enough to show that 
\be d_Z^S(\tilde{p_1},\tilde{p_2})\leq |(\tilde{x_1},0)-(y_1,H)| + ||y_1-y_2||_F+ |(y_2,H)-(\tilde{x_2},0)|.\ee

Notice
\begin{eqnarray*}
d_Z^S(\tilde{p_1},\tilde{p_2}) &\leq& d_Z^S(\tilde{p_1},\underline{q_1})+ d_Z^S(\underline{q_1},\underline{q_2})+d_Z^S(\underline{q_2},\tilde{p_2}) \quad \mbox{(Triangle Inequality)}\\
&\leq& |(\tilde{x_1},0)-(y_1,0)|+ d_Z^S(\underline{q_1},\underline{q_2})+|(y_2,0)-(\tilde{x_2},0)| \\
&\leq& |(\tilde{x_1},0)-(y_1,0)|+ ||y_1-y_2||_F+  \frac{K}{R} + |(y_2,0)-(\tilde{x_2},0)|\\ 
&\leq& |(\tilde{x_1},0)-(y_1,H)| +||y_1-y_2||_F+  \frac{K}{R} + |(y_2,0)-(\tilde{x_2},0)|\\
&\leq& |(\tilde{x_1},0)-(y_1,H)| +||y_1-y_2||_F+  |(y_2,H)-(\tilde{x_2},0)|
\end{eqnarray*}
where we use the inequality in theorem \ref{thm-SWIF=GH-R} and the last inequality is an application of \ref{H-choice}.

\subsection{Proof of GH convergence as $R\to\infty$}

%Hindy and Moshe and Sormani
\be
d_{GH}\left((\bar{B}^{R}_{r}(x_0), d_{R}), (F^{-1}([0,r]), d_F) \right) \to 0
\ee
We need only show the Hausdorff distances:
\be
d_{H}\left(\varphi_S(\bar{B}^R_{r}(x_0)), \varphi_H(F^{-1}([0,r])) \right) \le \delta_R
\ee
where $\delta_R \to 0$.

First we show
\be
\varphi_S(\bar{B}^R_{r}(x_0)) \subset T_\delta(\varphi_H(F^{-1}([0,r]))
\ee
where
\be\label{delta}
\delta= H + K/R
\ee
Taking any $x \in \bar{B}^R_{r}(x_0)$ we know from the set up that there
exists 
\be
v\in U^R_r(x_0)\subset \{v:\,\, ||v||_F < r +(K/R)\}\subset B_{2\lambda r}(0)\subset \mathbb{E}^N.
\ee
 such that $\pi_R(v)=x$.
Let
\be
w= rv/(r+K/R) \in F^{-1}[0,r] 
\ee
so we have
\be
||w||_F < ( r/(r+K/R) )(r +(K/R))=r.
\ee
%%Thus
%%\begin{eqnarray}
%%d_Z^R( \varphi_S(x), \varphi_F(w)) &=& d_Z^R(\pi_R(v,0) , \pi_R(w,H))\\
%&\le & d_Z^R(\pi_Z(v,0) , \pi_Z(w,0))+ d_Z^R(\pi_Z(w,0) , \pi_Z(w,H))\\
%&= & d_R(v,w)+ H  \\
%&\le& |v|\left(1-\frac{r}{r + K/R}\right) + H \\
%&\leq& 2\lambda r\left(\frac{K}{Rr + K}\right) + H = \delta 
%\end{eqnarray}

Thus
\begin{eqnarray}
d_Z^R( \varphi_S(x), \varphi_F(w)) &=& d_Z^R(\pi_Z(v,0) , \pi_Z(w,H))\\
&\le & d_Z^R(\pi_Z(v,0) , \pi_Z(v,H))+ d_Z^R(\pi_Z(v,H) , \pi_Z(w,H))\\
&\le & H + ||v - w||_F  \\
&\le& H + \left(1-\frac{r}{r + K/R}\right)||v||_F  \\
&\leq& H + \left(1-\frac{r}{r + K/R}\right)(r + K/R)\\
&=& H + K/R = \delta 
\end{eqnarray}

Next we show
\be
\varphi_H(F^{-1}([0,r]) \subset T_\delta(\varphi_S(\bar{B}_{rR}(x_0))): %%\text{typo with subscript $B_{rR}$. Should be B_{r}^R??}
\ee
Consider any $v \in F^{-1}([0,r]).$ 
Let
\be
w = \frac{(r- K/R)v}{r} \in F^{-1}[0,r] 
\ee
so we have
\be
||w||_F < \frac{(r- K/R)}{r}r = r - K/R.
\ee

By the definition of $K$, we have 
\be
|\bar{d_R}(w,0) - ||w||_F| \leq K/R
\ee
which implies that 
\be 
\bar{d_R}(w,0) \leq r.
\ee
Therefore, $x =\pi_R(w) \in \bar{B}_r^R(x_0).$
Thus
\begin{eqnarray}
d_Z^R( \varphi_S(x), \varphi_F(v)) &=& d_Z^R(\pi_Z(w,0) , \pi_Z(v,H))\\
&\le & d_Z^R(\pi_Z(w,0) , \pi_Z(w,H))+ d_Z^R(\pi_Z(w,H) , \pi_Z(v,H))\\
&\le&  H + ||v-w||_F \\
&=& H + \left(1 - \frac{(r - K/R)}{r}\right)||v||_F\\
&\leq& H + \left(1 - \frac{(r - K/R)}{r}\right)r \\
&=& H + K/R = \delta.
\end{eqnarray}
%%IS \delta' just H and second part of proof not necessary??

Taking $\delta_R=\delta$ we have 
\be
d_{H}\left(\varphi_S(\bar{B}^R_{r}(x_0)), \varphi_H(F^{-1}([0,r])) \right) \le \delta_R
\ee
and
\be
\lim_{R\to \infty} \delta_R= \lim_{R\to \infty} (H + K/R)
\ee
which converges to $0$ as $R\to \infty$.  
So
\be
d_{GH}\left((\bar{B}^{R}_{r}(x_0), d_{R}), (F^{-1}([0,r]), d_F) \right) \to 0
\ee

\subsection{Proof of SWIF Convergence as $R\to \infty$}

To prove
\be
d_{SWIF}\left((\bar{B}^{R}_{r}(x_0), d_{R}, T_R), (F^{-1}([0,r]), d_F, T_F) \right)\le M_R \to 0
\ee
we use the same distance preserving maps, $\varphi_H$ and $\varphi_S$
into the same $(Z_R, d^R_Z)$.  

So we need only prove there exists integral currents $A$ and $B$ such that
\be
\partial B + A =  \varphi_{S\#}T_X - \varphi_{H\#} T_F.
\ee
with
\be
\mass(B) + \mass(A) \le M_R
\ee
so that
\be
d_F^Z( \varphi_{S\#}T_X, \varphi_{H\#} T_F) \le \mass(B) + \mass(A)\le M_R
\ee
and we will have our claim.

We begin by defining the current $B$ which is depicted as the
lightly shaded diamond prism in Figure~\ref{fig:SWIF-Z}.

Recall that by the definition of $Z_R$ our smocking map
\be
\pi_Z: \, B_{2\lambda r}(0) \times [0,H]\subset {\mathbb{E}}^{N+1} \to Z_R
\ee
is surjective and by the definition of $d_Z^R$ and $d_Z^S$
is distance non-increasing
\be
d_Z^R(\pi_Z(v_1,v_2)) \le d_Z^R(\pi_Z(v_1,v_2))\le |v_1-v_2|.
\ee
Thus it is Lipschitz with $\dil(\pi_Z)\le 1$.   We apply it to define $B$ using a single chart
\be
B=\pi_{Z\#} \lbrack F^{-1}([0,r]) \times [0,H] \rbrack.
\ee
Then
\begin{eqnarray}
\mass(B) &\le& \dil(\pi_Z)^{N+1} \vol (F^{-1}([0,r]) \times [0,H])\\
&\le& H \vol (F^{-1}([0,r]) \,\,\le H (\dil(F^{-1})^N \omega_N r^N.
\end{eqnarray}
The boundary of $B$ is
\be
\partial B \,\,=\,\, S_{top}\,\,+\,\,S_{bottom}\,\,+\,\,S_{around}
\ee
where
\begin{eqnarray}
S_{top}&=&\pi_{Z\#} \lbrack F^{-1}([0,r]) \times \{H\} \rbrack \\
S_{bottom}&= &-\,\pi_{Z\#} \lbrack F^{-1}([0,r]) \times \{0\} \rbrack \\
S_{around}&=&\pi_{Z\#} \lbrack F^{-1}(\{r\}) \times [0,H] \rbrack. 
\end{eqnarray}
Note that
\be
S_{bottom}= \varphi_{H\#} T_F \textrm{ where } T_F=\lbrack F^{-1}([0,r]) \rbrack
\ee
is the integral current structure on $F^([0,r)$ with weight 1.

If $T_X$ is the standard integral current structure on 
the rescaled closed ball $\bar{B}_{r}(x_0)$ is our smocked metric space, then
\begin{eqnarray}
\varphi_{S\#}T_X&=&\varphi_{S\#}\pi_{X\#} \lbrack \pi_X^{-1}(\bar{B}_{r}(x_0)\rbrack \\
&=& \pi_{Z\#} \lbrack  \pi_X^{-1}(\bar{B}_r(x_0)\times \{0\} \rbrack .
\end{eqnarray}
Since this is not $S_{top}$ we set
\be
A_{top} \,\,=\,\, \varphi_{S\#}T_X \,\,-\,\, S_{top} = \pi_{Z\#} \lbrack U_1 \times \{0\} \rbrack
- \pi_{Z\#} \lbrack U_2 \times \{0\} \rbrack
\ee
where
\be
U_1= \{v \in {\mathbb{E}}^N: F(v)>r \textrm{ and } \bar{d}^R_X(v,0)/R<r\} 
\ee
and
\be
U_2= \{v \in {\mathbb{E}}^N: F(v)<r \textrm{ and } \bar{d}^R_X(v,0)/R>r\}. 
\ee
Since $\dil(\pi_Z)\le 1$ we can estimate the mass
\be
\mass(A_{top}) \le \vol(U_1) + \vol(U_2).
\ee
Setting $A= A_{top} - S_{around}$ (which is a disjoint difference) we have
\be
\mass(A)\le\vol(U_1) + \vol(U_2) + H \vol(F^{-1}(r))
\ee
which gives us
\be
\partial B + A =  \varphi_{S\#}T_X - \varphi_{H\#} T_F.
\ee
Now take
\begin{eqnarray}
M_R&=&\mass(A) + \mass(B)\\ 
&\le& \vol(U_1) + \vol(U_2) 
 + H \vol(F^{-1}(r))+ H (\dil(F^{-1})^N \omega_N r^N
\end{eqnarray}
We need only show this converges to $0$ for fixed $r$ and $F$ as $R\to \infty$. 

Recall that by (\ref{choice-H})
\be
\lim_{R\to \infty}  H= \lim_{R\to \infty} \sqrt{ 4rK/R + (K/R)^2}=0
 \ee
 so the first and last terms which are constants multiplied by $H$
 converge to $0$.
 
 Finally we estimate $\vol(U_i)$ using the hypotheses to show that
 \be
U_1\subset F^{-1}(r, r+(K/R)) \textrm{ and }
U_2\subset F^{-1}(r-(K/R), r). 
\ee
So that
\begin{eqnarray}
\vol(U_1)+\vol(U_2) &\le& \vol\left(F^{-1}(r-(K/R), r+(K/R))\right)\\
 &\le& 2\,(K/R)\,\dil(F^{-1}) \,\vol(F^{-1}(r)) \to 0. 
\end{eqnarray}
Thus we have proved intrinsic flat convergence.

This completes the proof of Theorem~\ref{thm-SWIF=GH-R}.   $\hspace{5.5cm} \square$

\section{The SWIF Tangent Cones of $X_+$, $X_\square$ and $X_T$}

We can now apply Theorem~\ref{thm-SWIF=GH-R} to find the
SWIF tangent cones of $X_+$, $X_\square$, and $X_T$.      
Recall that in Remark~\ref{all-four-nice} we explained that these
are all nice smocked spaces.   While it is not easy to find the norm, $||\cdot||_F$
for a given smocked metric space, and it requires a length double inductive proof to
estimate $\bar{d}_+(x, x')$, this has already been done for these three spaces.

\begin{ex}
The SWIF and GH tangent cone at $\infty$ of $X_+$ is a taxicab plane $({\mathbb{R}}^2, ||\cdot||_{F_+})$
where 
\be
F_+(x)= (|x_1|+|x_2|)/3
\ee
because it was proven by Shanell George, Vishnu Rendla, and Hindy Drillick in \cite{Smocked} that
\be
|\,\bar{d}_+(x, x')\,- \,[F_+(x)-F_+(x')] \,| \,\le \, K \qquad \forall x,x' \in {\mathbb{E}}^N.
\ee
\end{ex}

\begin{ex}
The SWIF and GH tangent cone at $\infty$ of $X_T$ is also isometric to a taxicab plane, 
$({\mathbb{R}}^2, ||\cdot||_{F_T})$,
where 
\be
F_T(x)= (|x_1|+|x_2|)/2
\ee
because it was proven by Kazaras, Dinowitz, and Afrifa in \cite{Smocked} that
\be
|\,\bar{d}_T(x, x')\,- \,[F_T(x)-F_T(x')] \,| \,\le \, K \qquad \forall x,x' \in {\mathbb{E}}^N.
\ee
\end{ex}

\begin{ex}
The SWIF and GH tangent cone at $\infty$ of $X_\square$ is a normed space whose unit ball is
an octagon, $({\mathbb{R}}^2, ||\cdot||_{F_\square})$,
where 
\be
F_\square(x)= 2\sqrt{2}\min\{|x_1|/3, |x_2|/3\} + 2 | |x_1|/3 - |x_2|/3 |
\ee
because it was proven by Huynh, Minichiello, and Hepburn in \cite{Smocked} that
\be
|\,\bar{d}_+(x, x')\,- \,[F_\square(x)-F_\square(x')] \,| \,\le \, K \qquad \forall x,x' \in {\mathbb{E}}^N.
\ee
\end{ex}

\begin{ex}
Drillick and Mujo discovered a smocked space whose tangent cone at 
infinity is a normed space whose unit ball is the convex hull of a countable collection of points located where lines of integer slope cross the ellipse $(x/2)^2 + y^2=1$.   Their smocked space, $X_=$, is defined by a periodic lattice of unit length horizontal line segments with
left endpoints located at $\{(2i, j):\, i,j \in {\mathbb{Z}}\}$ as in \cite{Smocked}.   
We believe this tangent cone is
also a SWIF limit but do not have access to their estimates to prove this.  We cite their presentation \cite{Smocked-Equals}.
\end{ex}

\begin{ex}
Other sequences of smocked metric spaces are studied in work of Antonetti, Farahzad, and Yamin
\cite{Smocked-H}.   We expect their limits should also be SWIF limits but leave the verification to them.
\end{ex}

A list of open problems on the various limits of smocked metric spaces has been included at the end of \cite{Smocked} some of which are at the level of undergraduates while others are more advanced and might be tackled by doctoral students.   Any work on the construction of three dimensional smocked metric spaces would be of significant interest as such spaces could be applied to better understand the SWIF limits of sequences of manifolds with almost nonnegative scalar curvature  \cite{Sormani-Scalar-21}.

\end{document}